\documentclass{gtart}

\def\ifplaintex{\expandafter\ifx\csname documentclass\endcsname\relax}

\def\gtp{{\mathsurround=0pt\it $\cal G\mskip-2mu$eometry \&\ 
$\cal T\!\!$opology $\cal P\!$ublications}}  

\def\recd{{\small Received:\qua\receiveddate\ifx\reviseddate\relax
\else\qquad Revised:\qua\reviseddate\fi\par}} 

\def\lognumber#1{\def\thelognumber{#1}}
\def\volumenumber#1{\def\thevolumenumber{#1}}
\def\volumeyear#1{\def\thevolumeyear{#1}}
\def\papernumber#1{\def\thepapernumber{#1}}
\def\pagenumbers#1#2{\def\startpage{#1}\def\finishpage{#2}}
\def\published#1{\def\publishdate{#1}}

\def\received#1{\def\receiveddate{#1}}
\def\revised#1{\def\reviseddate{#1}}
\def\accepted#1{\def\accepteddate{#1}}
\def\asciititle#1{\def\theasciititle{#1}}

\long\def\asciiabstract#1{\long\def\theasciiabstract{#1}}

\let\\\par\let\thelognumber\relax\let\thevolumenumber\relax
\let\thepapernumber\relax\let\thevolumeyear\relax\let\startpage\relax
\let\finishpage\relax\let\publishdate\relax\let\receiveddate\relax
\let\reviseddate\relax\let\accepteddate\relax\let\theasciititle\relax
\let\theasciiauthors\relax
\let\theasciiabstract\relax

\let\theasciiemail\relax

\ifplaintex
\font\logobig=cmssbx10 scaled 3836
\font\logomed=cmssbx10 scaled 2557
\else
\font\logobig=cmssbx10 scaled 4200
\font\logomed=cmssbx10 scaled 2800
\fi

\long\def\makeagttitle{   
\count0=\startpage
\agt\hfill      
\hbox to 45truept{\vbox to 0pt{\vglue -13truept{\logomed A\kern -.37em{\logobig 
T}\kern -.38em G}\vss}\hss}
\break
{\small Volume \thevolumenumber\ (\thevolumeyear)
\startpage--\finishpage\nl
Published: \publishdate}

\vglue .25truein

{\parskip=0pt\leftskip 0pt plus
1fil\def\\{\par\smallskip}{\Large\bf\thetitle}\par\medskip} \vglue
0.05truein

{\parskip=0pt\leftskip 0pt plus 1fil\def\\{\par}{\sc\theauthors}
\par\medskip}%
 
\vglue 0.03truein 

{\small\leftskip 25truept\rightskip 25truept{\bf Abstract}\stdspace\theabstract

{\bf AMS Classification}\stdspace\theprimaryclass
\ifx\thesecondaryclass\relax\else; \thesecondaryclass\fi\par
{\bf Keywords}\stdspace \thekeywords\par}\vglue 7truept

}   

\ifplaintex
\font\phead=cmsl9 scaled 950
\font\pnum=cmbx10 scaled 913
\font\pfoot=cmsl9 scaled 950
\headline{\vbox to 0pt{\vskip -4.5mm\line{\small\phead\ifnum
\count0=\startpage ISSN 1472-2739 (on-line) 1472-2747 (printed)
\hfill {\pnum\folio}\else\ifodd\count0\def\\{ }%
\ifx\theshorttitle\relax\thetitle\else\theshorttitle\fi\hfill{\pnum\folio}
\else\def\\{ and }{\pnum\folio}\hfill\ifx\theshortauthors\relax\theauthors
\else\theshortauthors\fi\fi\fi}\vss}}
\footline{\vbox to 0pt{\vglue 0mm\line{\small\pfoot\ifnum\count0=\startpage
\copyright\ \gtp\hfill\else
\agt, Volume \thevolumenumber\ (\thevolumeyear)\hfill\fi}\vss}}
\else
\font=cmsl9 scaled 1050
\font\lnum=cmbx10 
\font=cmsl9 scaled 1050
\makeatletter
\def\@oddhead{{\small\lhead\ifnum\count0=\startpage ISSN 1472-2739 
(on-line) 1472-2747 (printed)\hfill {\lnum\number\count0}\else\ifodd\count0
\def\\{ }\ifx\theshorttitle\relax \thetitle \else\theshorttitle\fi\hfill
{\lnum\number\count0}\else\def\\{ and }{\lnum\number\count0}
\hfill\ifx\theshortauthors\relax 
\theauthors\else\theshortauthors\fi\fi\fi}}\def\@evenhead{\@oddhead}
\def\@oddfoot{\small\lfoot\ifnum\count0=\startpage\copyright\ \gtp\hfill\else
\agt, Volume \thevolumenumber\ (\thevolumeyear)\hfill\fi}
\def\@evenfoot{\@oddfoot}
\makeatother
\fi
\let\maketitlepage\makeagttitle
\let\makeshorttitle\maketitlepage
\let\maketitle\maketitlepage

\newwrite\gtoutfile
\long\gdef\makeheadfile{  
{\def\\{, }\def\s{ }
\immediate\openout\gtoutfile head.xxx
\immediate\write\gtoutfile{Proxy-for: \ifx\theasciiauthors\relax
\theauthors\else\theasciiauthors\fi\s<\ifx\theasciiemail\relax\theemail\else\theasciiemail\fi>}
\immediate\write\gtoutfile{\noexpand\\}
\immediate\write\gtoutfile{Authors: \ifx\theasciiauthors\relax
\theauthors\else\theasciiauthors\fi}
{\def\\{ }\immediate\write\gtoutfile{Title: \ifx\theasciititle\relax
\thetitle\else\theasciititle\fi}}
\immediate\write\gtoutfile{Subj-class: GT or SG, GR etc}
\immediate\write\gtoutfile{MSC-class: \theprimaryclass\ifx\thesecondaryclass\relax\else, \thesecondaryclass\fi}
\immediate\write\gtoutfile{Journal-ref: Algebr. Geom. Topol. \thevolumenumber\s
(\thevolumeyear) \startpage-\finishpage}
\immediate\write\gtoutfile{Comments: Published by Algebraic and
Geometric Topology at}
\immediate\write\gtoutfile{\s\s\s  http://www.maths.warwick.ac.uk/agt/AGTVol\thevolumenumber/agt-\thevolumenumber-\thepapernumber.abs.html}
\immediate\write\gtoutfile{\noexpand\\}
\immediate\write\gtoutfile{}
\ifx\theasciiabstract\relax
\immediate\write\gtoutfile{\theabstract}\else
\immediate\write\gtoutfile{\theasciiabstract}\fi
\immediate\write\gtoutfile{}
\immediate\write\gtoutfile{\noexpand\\}
\immediate\write\gtoutfile{}
\immediate\closeout\gtoutfile}}  

\def\maketitlepage{\makeagttitle\makeheadfile}
\let\makeshorttitle\maketitlepage
\let\maketitle\maketitlepage

\lognumber{43}
\volumenumber{3}
\volumeyear{2003}
\papernumber{43}
\published{13 December 2003}
\pagenumbers{1225}{1256}
\received{26 February 2003}
\revised{8 December 2003}
\accepted{12 December 2003}

\usepackage{amssymb, amsmath}
\usepackage{amscd}
\usepackage{curves}
\usepackage{amsbsy}

\providecommand{\cal}{\mathcal}

\theoremstyle{plain} 
\newtheorem{theorem}{Theorem}[section] 
\newtheorem{lemma}[theorem]{Lemma}
\newtheorem{proposition}[theorem]{Proposition}
\newtheorem{corollary}[theorem]{Corollary}

\newtheorem{question}[theorem]{Question}

\newtheorem{conjecture}[theorem]{Conjecture}

\theoremstyle{definition} 

\theoremstyle{remark} 
\newtheorem{remark}[theorem]{Remark}
\newtheorem{example}[theorem]{Example}
\swapnumbers
\newtheorem{piece}[theorem]{}

\newenvironment{Proof}[1]{\begin{proof}[Proof of #1]}{\end{proof}} 

\newcommand{\refref}[1]{~\ref{#1} (p.~\pageref{#1})} 
\newcommand{\Secref}[1]{Sec\-tion \ref{#1}}
\newcommand{\Lemmaref}[1]{Lem\-ma \ref{#1}}
\newcommand{\Propref}[1]{Pro\-po\-si\-tion \ref{#1}} 
\newcommand{\Thmref}[1]{Theo\-rem \ref{#1}}
\newcommand{\Rkref}[1]{Re\-mark \ref{#1}}
\newcommand{\Corref}[1]{Co\-rol\-la\-ry \ref{#1}}
\newcommand{\Conjref}[1]{Con\-jec\-ture \ref{#1}}
\newcommand{\Figref}[1]{Fi\-gu\-re \ref{#1}} 
\newcommand{\Exref}[1]{Ex\-am\-ple \ref{#1}}

\newcommand{\thmref}[1]{Theo\-rem\refref{#1}}

\newcommand{\mraise}[2]{\raisebox{#1}{$#2$}}      

\newcommand{\sub}[1]{\raisebox{-2pt}{$\!_{#1}$} }  
\newcommand{\Sub}[1]{\raisebox{-2pt}{$\!_{\,#1}$} }
\newcommand{\adjustnabla}[1]{\sub{#1}} 
\newcommand{\goth}[1]{\mathfrak{#1}}              
\newcommand{\defemph}[1]{{\sffamily\slshape #1\/}}  

\renewcommand{\phi}{\varphi}
\renewcommand{\theta}{\vartheta}
\renewcommand{\epsilon}{\varepsilon} 

\newcommand{\ie}{{\it i.e.}~}
\newcommand{\eg}{{\it e.g.}~}
\newcommand{\st}{such that}
\newcommand{\Iff}{if and only if}

\newcommand{\LC}{Levi-Civit\`a}
\newcommand{\SW}{Sei\-berg--Witten}

\newcommand{\ac}{al\-most-com\-plex}
\newcommand{\riem}{Riemannian}

\newcommand{\str}{struc\-ture}

\newcommand{\aR}{\mathbb{R}}
\newcommand{\Zi}{\mathbb{Z}}

\newcommand{\Ce}{\mathbb{C}}

\newcommand{\Ha}{\mathbb{H}} 
\newcommand{\Proj}{\mathbb{P}}

\newcommand{\CP}{\Ce\Proj}
\newcommand{\Sph}[1]{{\mathbb S}^{#1}}  

\newcommand{\del}{\partial}           
\newcommand{\til}[1]{\widetilde{#1}}  
\renewcommand{\Bar}[1]{\overline{#1}} 
\newcommand{\rec}[1]{\tfrac{1}{#1}}   
\newcommand{\iso}{\approx}            
\newcommand{\maps}{\longmapsto}       
\newcommand{\inner}[1]{\langle #1\rangle} 
\newcommand{\Inner}[1]{\bigl\langle #1\bigr\rangle} 
\newcommand{\rest}[1]{|_{#1}}                          

\newcommand{\T}[1]{T_{#1}}             
\newcommand{\N}[1]{N_{#1}}             

\newcommand{\uaR}{\underline{\aR}}     
\newcommand{\Nabla}[1]{\nabla\adjustnabla{#1}}

\newcommand{\pinc}{pin$^{\!{\bf C}}$} 
\newcommand{\spinc}{s\pinc}         

\newcommand{\tsum}{\operatorname{\mraise{1pt}{\sum}}} 

\newcommand{\F}{{\mathcal F}}
\newcommand{\G}{{\mathcal G}}

\newcommand{\comment}[1]{%
	#1
	}

\begin{document}

\title{Existence of foliations\\on 4--manifolds}
\asciititle{Existence of foliations on 4-manifolds}
\author{Alexandru Scorpan}
\address{Department of Mathematics, University of 
Florida\\ 358 Little Hall, Gainesville, FL 32611--8105, USA}
\email{ascorpan@math.ufl.edu}

\url{www.math.ufl.edu/\char'176ascorpan}

\primaryclass{57R30}
\secondaryclass{57N13, 32Q60}
\keywords{Foliation, four-manifold, \ac}
\keywords{Foliation, four-manifold, almost-complex}

\begin{abstract}
We present existence results for certain singular $2$--dimensional
foliations on $4$--mani\-folds. The singularities can be chosen to be
simple, for example the same as those that appear in Lefschetz
pencils. There is a wealth of such creatures on most $4$--manifolds,
and they are rather flexible: in many cases, one can prescribe
surfaces to be transverse or be leaves of these foliations.

The purpose of this paper is to offer objects, hoping for a future
theory to be developed on them. For example, foliations that are taut
might offer genus bounds for embedded surfaces (Kronheimer's
conjecture).
\end{abstract}

\asciiabstract{We present existence results for certain singular 
2-dimensional foliations on 4-manifolds. The singularities can be
chosen to be simple, for example the same as those that appear in
Lefschetz pencils. There is a wealth of such creatures on most
4-manifolds, and they are rather flexible: in many cases, one can
prescribe surfaces to be transverse or be leaves of these foliations.

The purpose of this paper is to offer objects, hoping for a future
theory to be developed on them. For example, foliations that are taut
might offer genus bounds for embedded surfaces (Kronheimer's
conjecture).}

\maketitle

\section{Introduction}

Foliations play a very important role in the study of $3$--manifolds, but almost none so far in the study of $4$--manifolds.  There are hints, though, that they should play an important role here as well. For example, for $M^4=N^3\times\Sph{1}$, Kronheimer obtained genus bounds for embedded surfaces from certain taut foliations \cite{kronheimer}, which are sharper than the ones coming from \SW\ basic classes. He conjectured that such bounds might hold in general.

\subsection{Summary}

(In this paper, all foliations will be $2$--dimensional and oriented, all manifolds will be $4$--dimensional, closed and oriented; unless otherwise specified, of course.)

For a foliation $\F$ to exist on a manifold $M$, the tangent bundle must split $\T{M}=\T{\F}\oplus\N{\F}$. Since in general that does not happen, one must allow for singularities of $\F$. An important example is \cite{donaldson.lefschetz}:

\begin{example}[S.K.\ Donaldson 1999]
Let $J$ be an \ac\ \str\ on $M$ that admits a compatible symplectic structure (\ie $J$ admits a closed $2$--form $\omega$ \st\ $\omega(x,Jx)>0$ and $\omega(Jx,Jy)=\omega(x,y)$). Then $J$ can be deformed to an \ac\ \str\ $J'$ \st\ $M$ admits a Lefschetz pencil with $J'$--holomorphic fibers. 
\end{example}

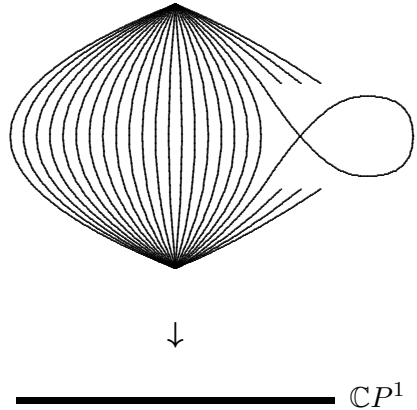
\begin{figure}[t]
\begin{center}
\setlength{\unitlength}{0.5pt}
\begin{picture}(300,320)(0,0)

\qbezier(150,300)(-100,200)(150,100)
\qbezier(150,300)(-80,200)(150,100)
\qbezier(150,300)(-60,200)(150,100)
\qbezier(150,300)(-40,200)(150,100)
\qbezier(150,300)(-20,200)(150,100)
\qbezier(150,300)(0,200)(150,100)
\qbezier(150,300)(20,200)(150,100)
\qbezier(150,300)(40,200)(150,100)
\qbezier(150,300)(60,200)(150,100)
\qbezier(150,300)(80,200)(150,100)
\qbezier(150,300)(100,200)(150,100)
\qbezier(150,300)(120,200)(150,100)
\qbezier(150,300)(140,200)(150,100)
\qbezier(150,300)(160,200)(150,100)
\qbezier(150,300)(180,200)(150,100)
\qbezier(150,300)(200,200)(150,100)
\qbezier(150,300)(220,200)(150,100)
\qbezier(150,300)(240,200)(150,100)
\qbezier(150,300)(260,200)(150,100)
\qbezier(150,300)(280,200)(150,100)

\qbezier(150,300)(185,280)(215,240)
\qbezier(150,100)(185,120)(215,160)

\qbezier(215,240)(260,170)(290,170)
\qbezier(215,160)(260,230)(290,230)

\qbezier(290,170)(330,167)(330,200)
\qbezier(290,230)(330,233)(330,200)

\qbezier(150,300)(190,280)(230,240)
\qbezier(150,100)(190,120)(230,160)

\qbezier(150,300)(195,280)(245,240)
\qbezier(150,100)(195,120)(245,160)

\qbezier(150,300)(200,280)(260,240)
\qbezier(150,100)(200,120)(260,160)

\put(150,50){\makebox(0,0){$\boldsymbol{\downarrow}$}}

\linethickness{2pt}
\put(30,0){\line(1,0){240}}
\put(280,-5){$\mathbb{C}P^1$}

\end{picture}
\end{center}
\caption{A Lefschetz pencil}
\label{fig-lefschetz}
\end{figure}

A \defemph{Lefschetz pencil} is a \emph{singular fibration} $M\to\CP^1$ with singularities modeled locally by 
\[ (z_1,z_2)\maps z_1/z_2  \qquad\text{or}\qquad (z_1,z_2)\maps z_1 z_2 \]
for suitable local complex coordinates (compatible with the orientation of $M$). Note that all fibers pass through all singularities of type $z_1/z_2$, see \Figref{fig-lefschetz}. 
The existence of a Lefschetz pencil is equivalent to the existence of a symplectic structure. See \cite[ch.~8]{gompf} for a survey.

\medskip

A main result of this paper (\Thmref{thm-exist}) is that, under mild homological conditions on $M$, \emph{any} \ac\ \str\ $J$ on $M$ can be deformed to a $J'$ that admits a \emph{singular foliation} with $J'$--holomorphic leaves, and with singularities of the same type as those appearing in a Lefschetz pencil. In fact, the singularities can be chosen to be all of $z_1/z_2$--type (``pencil'' singularities), and thus could be eliminated by blow-ups. Other singularities can also be chosen (or even just a single complicated singularity), see \Secref{sec-other.sing}. 

By allowing singularities with reversed orientations, this existence result can be generalized to \spinc-\str s (\Thmref{thm-exist2}). (S\pinc-\str s are more general than \ac\ \str s and always exist on any $4$--manifold.) 

Also, under certain natural conditions, given embedded surfaces can be arranged to be transversal to the foliations (\Thmref{thm-transversal}), or even to be leaves of the foliations (\Thmref{thm-leaf}).

The main tools used in proving our results are: Thurston's $h$--princi\-ple  for foliations with codimension $\geq 2$ (see \ref{thm-thurston} below), which takes care of integrating plane-fields and reduces the existence problem to bundle theory; and the Dold--Whitney Theorem characterizing bundles by their characteristic classes (see \ref{thm-dold.whitney} below). In a nutshell, we build a bundle with the same characteristic classes as $\T{M}$, we conclude it is $\T{M}$, we let Thurston integrate to a foliation.

This paper presents existence results for singular $2$--di\-men\-sio\-nal foliations on $4$--manifolds. We offer a wealth of objects to be used in a future theory, and try to stimulate interest in this area.

The paper is organized as follows: \Secref{sec-statements} contains the statements of most results of this paper, \Secref{sec-context} offers a quick survey of the context in foliation theory, \Secref{sec-proofs} presents the proofs of our results, while \Secref{sec-appendix} contains left-overs.

\subsection{Why bother?}

One hint that foliations on $4$--manifolds are worth studying (especially \emph{taut} foliations, see \Secref{sec-taut} for a discussion) are Kronheimer's results (see \Thmref{thm-kron} and \Conjref{conj} below). Taut foliations might offer \emph{minimal genus bounds} for embedded surfaces, see \Secref{sec-min.genus}.

In a slightly larger context, the relationship (if any) between foliations on $4$--manifolds and \emph{\SW\ theory} is worth elucidating.

Another question worth asking is: For what foliations is the induced \ac\ \str\ ``nice'' (\ie close to symplectic). One such problem asks for which foliations does the induced \ac\ \str\ have \emph{Gromov compactness} (\ie whether the space of $J$--holomorphic curves of a fixed genus and homology class is compact; see 
Question~\ref{q-kirby}). 

In general, one can hope that foliations will help better visualize, manipulate and understand \emph{\ac\ \str s}, maybe in a manner similar to the one in which open-book decompositions help understand contact \str s on $3$--manifolds (see also \Corref{cor-no.gromov}).

\subsection*{Acknowledgments}

We wish to thank Rob Kirby for his constant encouragement and wise advice (mathematical and otherwise).

\section{Statements}
\label{sec-statements}

(In this paper, Poincar\'e duality will be used blindly, submanifolds and the homology classes they represent will frequently be denoted with the same symbol, and top (co)ho\-mo\-logy classes will be paired with fundamental cycles without comment. For example, $\chi(M)-\tau\nu$ could be written more elegantly as $\chi(M)-(\tau\cup\nu)[M]$, while $c_1(L)=\Sigma$ is $c_1(L)=PD([\Sigma])$.)

First of all, notice that any (non-singular) foliation $\F$ on $M$ induces \ac\ \str s: Pick a \riem\ metric $g$, embed the normal bundle $\N{\F}$ in $\T{M}$, then define an \ac\ \str\ $J_\F$ to be the rotation by $\pi/2$ (respecting orientations) in both $\T{\F}$ and $\N{\F}$. It has the property that the leaves of $\F$ are $J_\F$--holomorphic. (In general, we will call an \ac\ \str\
$J$ \defemph{compatible} with a foliation $\F$ if $J$ makes the leaves of $\F$ be $J$--holomorphic.)
The first Chern class $c_1(J_\F)$ is well-defined independent of the choices made. We have:
\begin{gather*}
c_1(J_\F)=e(\T{\F})+e(\N{\F})\\
\chi(M)=e(\T{\F})\cdot e(\N{\F})
\end{gather*}
If the foliation $\F$ has singularities, then the second equality above fails, and the defect $\chi(M)-e(\T{\F})\cdot e(\N{\F})$ measures the number of singularities (or, for more general singularities, their complexity, see \Thmref{thm-exist.sing}).

\subsection{Main existence results}
Call a class $c\in H^{2}(M;\Zi)$ a \defemph{complex class} of the $4$--manifold $M$ if
\[ c\equiv w_2(M) \pmod{2} \qquad\text{and}\qquad 
  p_1(M)=c^2-2\chi(M) \]

An element $c\in H^2(M;\Zi)$ is a  complex class \Iff\ there is an \ac\ \str\ $J$ on $M$ \st\ $c_1(\T{M},J)=c$. 
One direction is elementary: If $J$ is such a \str, then $(\T{M},J)$ is a complex-plane bundle, and thus has $c_1(\T{M})\equiv w_2(\T{M})\pmod{2}$ and $p_1(\T{M})=c_1(\T{M})-2c_2(\T{M})$. The converse was proved in \cite{wu,HH} (and will appear here re-proved as part of \Corref{cor-cx.str}).

\begin{theorem}[Existence Theorem]
\label{thm-exist}
Let $c\in H^2(M;\Zi)$ be a complex class, and let $c=\tau+\nu$ be any splitting \st\ $\chi(M)-\tau\nu\geq 0$.
Choose any combination of $n=\chi(M)-\tau\nu$ singularities modeled on the levels of the complex functions $(z_1,z_2)\maps z_1/z_2$ or $(z_1,z_2)\maps z_1 z_2$.
Then there is a singular foliation $\F$ with $e(\T{\F})=\tau$,  $e(\N{\F})=\nu$, and with $n$ singularities as prescribed.
\end{theorem}

\begin{remark}
\label{rk-euler}
Due to the singularities, the bundles $\T{\F}$ and $\N{\F}$ are only defined on $M\setminus\{\text{singularities}\}$. Their Euler classes \emph{a priori} belong to $H^2(M\setminus\{\text{singularities}\};\,\Zi)$, but can be pulled-back to $H^2(M;\Zi)$, since the isolated singularities can be chosen to affect only the 4-skeleton of $M$, and thus to not influence $H^2$.
\end{remark}

\begin{remark}
Unlike a Lefschetz pencil, in general \emph{not all} leaves of the foliation pass through the $z_1/z_2$--singularities. 
See \Exref{ex-create.torus} for creating a leaf that touches no singularity.
\end{remark}

\subsection{Restrictions}

Finding a splitting $c=\tau+\nu$ with $\chi(M)-\tau\nu\geq 0$ is possible for most $4$--manifolds that admit \ac\ \str s. For example, if 
\[ \chi(M)\geq 0 \] 
(\eg for all \emph{simply-connected} $M$'s), then one can choose either one of $\tau$ or $\nu$ to be $0$, and conclude that such foliations exist. Or:

\begin{lemma}
If $b_2^+(M)>0$, then there are infinitely many splittings $c=\tau+\nu$ with $\chi(M)-\tau\nu\geq 0$ (and thus infinitely many homotopy types of foliations). 
\end{lemma}

\begin{proof}
If $b_2^+(M)>0$, there is a class $\alpha$ with $\alpha\cdot\alpha>0$. Choose $\tau=c-k\alpha$  and $\nu=k\alpha$ ($k\in\Zi$). 
Then $\chi(M)-\tau\nu=\chi(M)-kc\alpha+k^2\alpha^2$, and for $k$ big enough it will be positive.
\end{proof}

The main restriction to the existence of such foliations remains, of course, the existence of a complex class.
But \Thmref{thm-exist} can be generalized for the case when $c$ is merely an \emph{integral lift} of $w_2(M)$, see \Thmref{thm-exist2}. 
In that case, singularities are also modeled using local complex coordinates, but are allowed to be compatible either with the orientation of $M$ or with the \emph{opposite orientation}.

This is similar to the generalization of Lefschetz pencils to
\emph{achiral Lefschetz pencils}, see \cite[\S8.4]{gompf} and \Secref{sec-not.cx}.
As it happens, the only \emph{known} obstruction to the existence of an achiral Lefschetz pencil (\cite[8.4.12--13]{gompf}) is the \emph{only} obstruction to the existence of such an ``achiral'' singular foliation
(see \Secref{sec-not.cx} and \Propref{prop-achiral.obstruction}).

\subsection{Singularities}

The singularities of $\F$ are exactly the singularities that appear in a Lefschetz pencil. They can be chosen in either combination of types as long as their number is $n=\chi(M)-\tau\nu$. For example, there are always foliations with only $z_1/z_2$--singularities, that can thus be eliminated by blowing-up.
In fact, other choices of singularities are possible.

Namely, for any isolated singularity $p$ of a foliation that is compatible with an \ac\ \str\ we will define its \defemph{Hopf degree} $\deg p\geq 0$ (essentially a Hopf invariant of the tangent plane field above a small $3$--sphere around $p$; see \Secref{sec-other.sing}).
Then: 

\begin{theorem}
\label{thm-exist.sing}
Let $c\in H^2(M;\Zi)$ be a complex class, and let $c=\tau+\nu$ be any splitting \st\ $\chi(M)-\tau\nu\geq 0$.
Then, for any choice of (positive) singularities  $\{p_1,\ldots,p_k\}$ so that
\[ \tsum \deg p_i = \chi(M)-\tau\nu \]
there is a singular foliation $\F$ with $e(\T{\F})=\tau$, $e(\N{\F})=\nu$, and with the chosen singularities.
\end{theorem}

In analogy to the Poincar\'e--Hopf theorem on indexes of vector fields, a converse to the above is true:

\begin{proposition}
\label{prop-index}
For any singular foliation $\F$ on $M$ with isolated singularities $\{p_1,\ldots,p_k\}$ compatible with a local \ac\ \str, we have
\[ \chi(M) = \tsum \deg p_i + e(\T{\F})\cdot e(\N{\F})  \]
\end{proposition}

\subsection{Prescribing leaves and closed transversals}

Let $\F$ be a foliation. If $S$ is a closed transversal of $\F$, then we must have 
\begin{align*}
& e(\T{\F})\cdot S=e(\T{\F}\rest{S})=e(\N{S})=S\cdot S\\
& e(\N{\F})\cdot S=e(\N{\F}\rest{S})=e(\T{S})=\chi(S)
\end{align*}
These conditions are, in fact, sufficient:

\begin{theorem}[Closed transversal]
\label{thm-transversal}
Let $S$ be a closed connected surface. 
Let $c$ be a complex class with a splitting  $c=\tau+\nu$ \st\ $\chi(M)-\tau\nu\geq 0$. 
If
\[ \chi(S)=\nu\cdot S \qquad\qquad 
  S\cdot S=\tau\cdot S \]
then there is a singular foliation $\F$ with $e(\T{\F})=\tau$, $e(\N{\F})=\nu$, and having $S$ as a closed transversal.
\end{theorem}

If, on the other hand, $S$ is a closed leaf of $\F$, then we have 
\begin{align*}
& e(\T{\F})\cdot S=e(\T{\F}\rest{S})=e(\T{S})=\chi(S)\\
& e(\N{\F})\cdot S=e(\N{\F}\rest{S})=e(\N{S})=S\cdot S
\end{align*}
Conversely:

\begin{theorem}[Closed leaf]
\label{thm-leaf}
Let $S$ be a closed connected surface with 
$S\cdot S\geq 0$. 
Let $c$ be a complex class with a splitting $c=\tau+\nu$ \st\  $\chi(M)-\tau\nu\geq S\cdot S$.
If
\[ \chi(S)=\tau\cdot S \qquad\qquad 
  S\cdot S=\nu\cdot S \]
then there is a singular foliation $\F$ with $e(\T{\F})=\tau$, $e(\N{\F})=\nu$, and having $S$ as a closed leaf.
(The number of singularities along $S$ is $S\cdot S$.)
\end{theorem}

(Surfaces with $S\cdot S<0$ could be made leaves of achiral singular foliations, using singularities with reversed orientations, see \ref{prop-achiral.leaf}.)

An immediate consequence of the above is:

\begin{corollary}[Trivial tori]
\label{cor-torus}
A homologically-trivial torus can always be made a leaf or a transversal of a foliation.
\end{corollary}

Such flexibility is a strong suggestion that more rigidity is needed in order to actually catch any of the topology of $M$ with the aid of foliations. Requiring foliations to be taut seems a natural suggestion. (Compare with \Exref{ex-create.torus}.)

The conditions $\chi(S)=\tau\cdot S$ and $S\cdot S=\nu\cdot S$ from \ref{thm-leaf} add to $\chi(S)+S\cdot S=c\cdot S$. The conditions $\chi(S)=\nu\cdot S$ and $S\cdot S=\tau\cdot S$ from \ref{thm-transversal} also add to $\chi(S)+S\cdot S=c\cdot S$. For good choices of $\tau$ and $\nu$, that is sufficient:

\begin{corollary}[Adjunct surfaces]
\label{cor-adjunct}
Let $c$ be a complex class, and let $S$ be a closed connected surface such that 
\[ \chi(S)+S\cdot S=c\cdot S \]
If $\chi(M)-\chi(S)\geq 0$, then there is a singular foliation $\F_1$ with $e(\T{\F_1})=S$ that has $S$ as a closed transversal.\\
If further $\chi(M)-\chi(S)\geq S\cdot S\geq 0$, then there is also a singular foliation $\F_2$ with $e(\N{\F_2})=S$ that has $S$ as a leaf (with $S\cdot S$ singularities on it).
\end{corollary}

\begin{proof}
For $\F_1$, pick $\tau=S$ and $\nu=c-S$. For $\F_2$, pick $\tau=c-S$ and $\nu=S$. (In both cases $\tau\nu=\chi(S)$, and so $\chi(M)-\tau\nu\geq 0$.) Apply \ref{thm-leaf} or \ref{thm-transversal}.
\end{proof}

As a consequence of (the proof of) \ref{cor-adjunct}, we can also re-prove the following \cite{bohr}:

\begin{proposition}[C.~Bohr 2000]
\label{prop-bohr}
Let $S$ be an embedded closed connected surface and $c$ a complex class.
Then there is an \ac\ \str\ $J$ with $c_1(J)=c$ \st\ $S$ is $J$--holomorphic
\Iff\ $\chi(S)+S\cdot S=c\cdot S$.
\end{proposition}

\begin{proof}
The positivity condition $\chi(M)-\chi(S)\geq 0$ is only needed for integrating the singularities of the foliations, and thus it can be ignored here. We have a (singular) plane field $\T{\F}$ that is transverse to $S$, and a (singular) plane field $\N{\F}$ that can be arranged to be tangent to $S$. These plane fields induce an \ac\ \str\ $J_\F$ that leaves $\T{S}$ invariant.
\end{proof}

\comment{
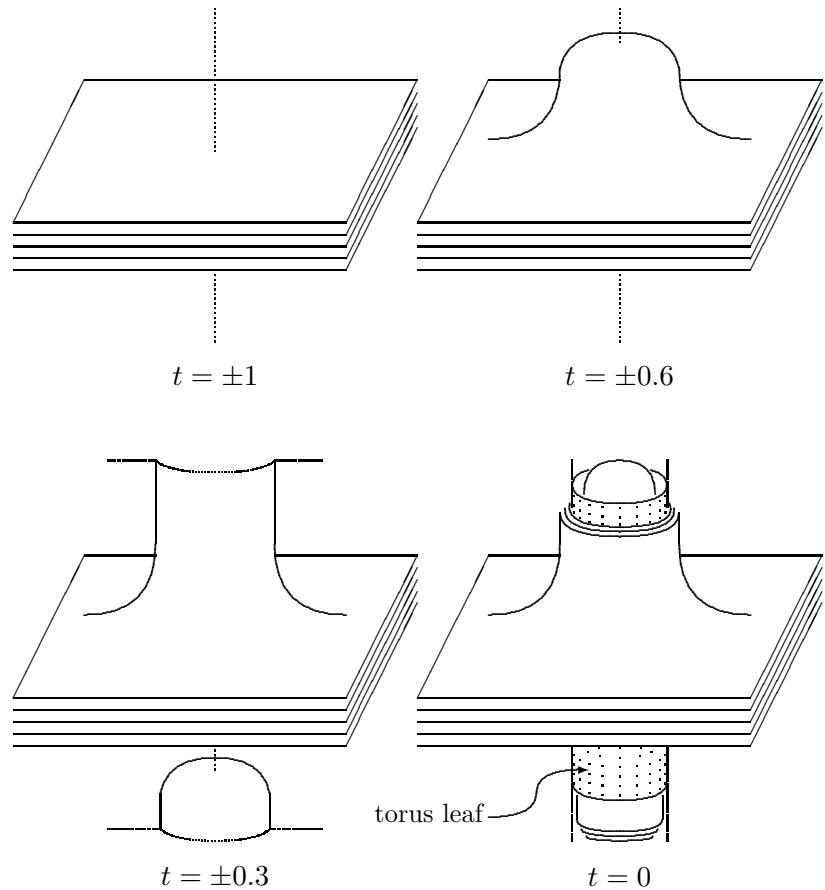
\begin{figure}[htb] 
\begin{center}
\setlength{\unitlength}{0.9pt}
\begin{picture}(340,370)(0,40)

\put(0,300){\line(1,0){140}}
\put(0,305){\line(1,0){140}}
\put(0,310){\line(1,0){140}}
\put(0,315){\line(1,0){140}}
\put(0,320){\line(1,0){140}}
\put(140,300){\line(1,2){30}}
\put(140,305){\line(1,2){30}}
\put(140,310){\line(1,2){30}}
\put(140,315){\line(1,2){30}}
\put(140,320){\line(1,2){30}}
\put(0,320){\line(1,2){30}}

\put(30,380){\line(1,0){140}}

\qbezier[30](85,350)(85,380)(85,410) 
\qbezier[15](85,300)(85,285)(85,270)

\put(85,255){\makebox(0,0){$t=\pm 1$}} 

\put(170,300){\line(1,0){140}} 
\put(170,305){\line(1,0){140}}
\put(170,310){\line(1,0){140}}
\put(170,315){\line(1,0){140}}
\put(170,320){\line(1,0){140}}
\put(310,300){\line(1,2){30}}
\put(310,305){\line(1,2){30}}
\put(310,310){\line(1,2){30}}
\put(310,315){\line(1,2){30}}
\put(310,320){\line(1,2){30}}
\put(170,320){\line(1,2){30}}

\put(200,380){\line(1,0){30}} 
\put(340,380){\line(-1,0){60}}
\qbezier(200,355)(230,355)(230,385)
\qbezier(310,355)(280,355)(280,385)

\qbezier(230,385)(233,400)(255,400) 
\qbezier(280,385)(277,400)(255,400)

\qbezier[7](255,396)(255,403)(255,410) 
\qbezier[15](255,300)(255,285)(255,270)

\put(255,255){\makebox(0,0){$t=\pm 0.6$}}       

\put(0,100){\line(1,0){140}} 
\put(0,105){\line(1,0){140}}
\put(0,110){\line(1,0){140}}
\put(0,115){\line(1,0){140}}
\put(0,120){\line(1,0){140}}
\put(140,100){\line(1,2){30}}
\put(140,105){\line(1,2){30}}
\put(140,110){\line(1,2){30}}
\put(140,115){\line(1,2){30}}
\put(140,120){\line(1,2){30}}
\put(0,120){\line(1,2){30}}

\put(30,180){\line(1,0){30}} 
\put(170,180){\line(-1,0){60}}
\qbezier(30,155)(60,155)(60,185)
\qbezier(140,155)(110,155)(110,185)

\put(60,185){\line(0,1){35}} 
\put(110,185){\line(0,1){35}}

\qbezier[20](40,220)(50,220)(60,220) 
\qbezier[20](110,220)(120,220)(130,220)
\qbezier[27](60,220)(65,215)(85,215)
\qbezier[27](85,215)(105,215)(110,220)

\put(62,65){\line(0,1){15}} 
\put(108,65){\line(0,1){15}}
\qbezier(62,80)(65,95)(85,95)
\qbezier(108,80)(105,95)(85,95)

\qbezier[22](40,65)(50,65)(62,65) 
\qbezier[22](108,65)(120,65)(130,65)
\qbezier[24](62,65)(67,60)(85,60)
\qbezier[24](85,60)(103,60)(108,65)

\qbezier[5](85,90)(85,95)(85,100) 

\put(85,45){\makebox(0,0){$t=\pm 0.3$}}    

\put(170,100){\line(1,0){140}} 
\put(170,105){\line(1,0){140}}
\put(170,110){\line(1,0){140}}
\put(170,115){\line(1,0){140}}
\put(170,120){\line(1,0){140}}
\put(310,100){\line(1,2){30}}
\put(310,105){\line(1,2){30}}
\put(310,110){\line(1,2){30}}
\put(310,115){\line(1,2){30}}
\put(310,120){\line(1,2){30}}
\put(170,120){\line(1,2){30}}

\put(200,180){\line(1,0){30}} 
\put(340,180){\line(-1,0){60}}
\qbezier(200,155)(230,155)(230,185)
\qbezier(310,155)(280,155)(280,185)

\put(230,185){\line(0,1){13}} 
\put(280,185){\line(0,1){13}}
\qbezier(230,198)(230,188)(255,188)
\qbezier(280,198)(280,188)(255,188)


\qbezier(232,200)(231,190)(255,190) 
\qbezier(278,200)(279,190)(255,190)

\qbezier(234,202)(232,192)(255,192) 
\qbezier(276,202)(280,192)(255,192)

\put(235,200){\line(0,1){10}} 
\put(275,200){\line(0,1){10}}
\qbezier(235,210)(235,202)(255,202)
\qbezier(275,210)(275,202)(255,202)
\qbezier(235,210)(236,215)(242,216)
\qbezier(275,210)(274,215)(268,216)

\qbezier[5](235,208)(235,200)(255,200) 
\qbezier[5](275,208)(275,200)(255,200)
\qbezier[5](235,204)(235,196)(255,196)
\qbezier[5](275,204)(275,196)(255,196)
\qbezier[5](235,201)(235,193)(255,193)
\qbezier[5](275,201)(275,193)(255,193)

\qbezier[10](235,210)(235,215)(235,220) 
\qbezier[10](275,210)(275,215)(275,220)

\qbezier(240,206)(240,220)(255,220) 
\qbezier(270,206)(270,220)(255,220)

\put(235,100){\line(0,-1){15}} 
\put(275,100){\line(0,-1){15}}
\qbezier(235,85)(235,77)(255,77)
\qbezier(275,85)(275,77)(255,77)

\qbezier[5](235,87)(235,79)(255,79) 
\qbezier[5](275,87)(275,79)(255,79)
\qbezier[5](235,91)(235,83)(255,83)
\qbezier[5](275,91)(275,83)(255,83)
\qbezier[5](235,95)(235,87)(255,87)
\qbezier[5](275,95)(275,87)(255,87)
\qbezier[5](235,99)(235,91)(255,91)
\qbezier[5](275,99)(275,91)(255,91)
\qbezier[3](240,100)(240,95)(255,95)
\qbezier[3](270,100)(270,95)(255,95)
\qbezier[2](250,100)(240,99)(255,99)
\qbezier[2](260,100)(270,99)(255,99)

\qbezier[25](235,85)(235,72.5)(235,60) 
\qbezier[25](275,85)(275,72.5)(275,60)

\put(237,69){\line(0,1){10}} 
\put(273,69){\line(0,1){10}}
\qbezier(237,69)(237,64)(255,64)
\qbezier(273,69)(273,64)(255,64)

\qbezier(239,65)(239,62)(255,62) 
\qbezier(271,65)(271,62)(255,62)

\qbezier(241,62)(241,60)(255,60) 
\qbezier(269,62)(269,60)(255,60)

\put(255,45){\makebox(0,0){$t=0$}} 
\put(242,90){\vector(1,0){0}}
\qbezier(239,90)(220,90)(220,80)
\qbezier(220,80)(220,70)(200,70)
\put(152,68){{\small torus leaf}}

\end{picture}
\end{center}
\caption{Creating a torus leaf (3D movie)}
\label{fig-create.torus1}
\end{figure}
}
\comment{
\begin{figure}[htb] 
\begin{center}
\setlength{\unitlength}{0.9pt}
\begin{picture}(300,275)(0,25)

\put(0,200){\line(1,0){140}}
\put(0,210){\line(1,0){140}}
\put(0,220){\line(1,0){140}}
\put(0,230){\line(1,0){140}}
\put(0,240){\line(1,0){140}}
\put(0,250){\line(1,0){140}}
\put(0,260){\line(1,0){140}}
\put(0,270){\line(1,0){140}}
\put(0,280){\line(1,0){140}}
\put(0,290){\line(1,0){140}}

\qbezier[55](70,190)(70,255)(70,300)

\put(52,180){$t=\pm 1$}

\put(160,200){\line(1,0){40}}
\put(160,210){\line(1,0){40}}
\put(160,220){\line(1,0){40}}
\put(160,230){\line(1,0){40}}
\put(160,240){\line(1,0){40}}
\put(160,250){\line(1,0){40}}
\put(160,260){\line(1,0){40}}
\put(160,270){\line(1,0){40}}
\put(160,280){\line(1,0){40}}
\put(160,290){\line(1,0){40}}

\put(300,200){\line(-1,0){40}}
\put(300,210){\line(-1,0){40}}
\put(300,220){\line(-1,0){40}}
\put(300,230){\line(-1,0){40}}
\put(300,240){\line(-1,0){40}}
\put(300,250){\line(-1,0){40}}
\put(300,260){\line(-1,0){40}}
\put(300,270){\line(-1,0){40}}
\put(300,280){\line(-1,0){40}}
\put(300,290){\line(-1,0){40}}

\qbezier(200,270)(210,270)(215,290)
\qbezier(245,290)(250,270)(260,270)

\qbezier(200,260)(210,260)(215,280)
\qbezier(245,280)(250,260)(260,260)

\qbezier(200,250)(210,250)(215,270)
\qbezier(215,270)(220,290)(230,290)
\qbezier(230,290)(240,290)(245,270)
\qbezier(245,270)(250,250)(260,250)

\qbezier(200,240)(210,240)(215,260)
\qbezier(215,260)(220,280)(230,280)
\qbezier(230,280)(240,280)(245,260)
\qbezier(245,260)(250,240)(260,240)

\qbezier(200,230)(210,230)(215,250)
\qbezier(215,250)(220,270)(230,270)
\qbezier(230,270)(240,270)(245,250)
\qbezier(245,250)(250,230)(260,230)

\qbezier(200,220)(210,220)(215,240)
\qbezier(215,240)(220,260)(230,260)
\qbezier(230,260)(240,260)(245,240)
\qbezier(245,240)(250,220)(260,220)

\qbezier(200,210)(210,210)(215,230)
\qbezier(215,230)(220,250)(230,250)
\qbezier(230,250)(240,250)(245,230)
\qbezier(245,230)(250,210)(260,210)

\qbezier(200,200)(210,200)(215,220)
\qbezier(215,220)(220,240)(230,240)
\qbezier(230,240)(240,240)(245,220)
\qbezier(245,220)(250,200)(260,200)

\qbezier(215,210)(220,230)(230,230)
\qbezier(230,230)(240,230)(245,210)

\qbezier(215,200)(220,220)(230,220)
\qbezier(230,220)(240,220)(245,200)

\qbezier[55](230,190)(230,255)(230,300)

\put(208,180){$t=\pm 0.6$}

\put(0,50){\line(1,0){40}}
\put(0,60){\line(1,0){40}}
\put(0,70){\line(1,0){40}}
\put(0,80){\line(1,0){40}}
\put(0,90){\line(1,0){40}}
\put(0,100){\line(1,0){40}}
\put(0,110){\line(1,0){40}}
\put(0,120){\line(1,0){40}}
\put(0,130){\line(1,0){40}}
\put(0,140){\line(1,0){40}}

\put(140,50){\line(-1,0){40}}
\put(140,60){\line(-1,0){40}}
\put(140,70){\line(-1,0){40}}
\put(140,80){\line(-1,0){40}}
\put(140,90){\line(-1,0){40}}
\put(140,100){\line(-1,0){40}}
\put(140,110){\line(-1,0){40}}
\put(140,120){\line(-1,0){40}}
\put(140,130){\line(-1,0){40}}
\put(140,140){\line(-1,0){40}}

\qbezier(40,50)(55,50)(55,140)
\qbezier(40,60)(54,60)(54,140)
\qbezier(40,70)(53,70)(53,140)
\qbezier(40,80)(52,80)(52,140)
\qbezier(40,90)(51,90)(51,140)
\qbezier(40,100)(50,100)(50,140)
\qbezier(40,110)(49,110)(49,140)
 
\qbezier(100,50)(85,50)(85,140)
\qbezier(100,60)(86,60)(86,140)
\qbezier(100,70)(87,70)(87,140)
\qbezier(100,80)(88,80)(88,140)
\qbezier(100,90)(89,90)(89,140)
\qbezier(100,100)(90,100)(90,140)
\qbezier(100,110)(91,110)(91,140)

\qbezier(54,50)(54,140)(70,140)
\qbezier(55,50)(55,130)(70,130)
\qbezier(56,50)(56,120)(70,120)
\qbezier(57,50)(57,110)(70,110)
\qbezier(58,50)(58,100)(70,100)
\qbezier(59,50)(59,90)(70,90)
\qbezier(60,50)(60,80)(70,80)

\qbezier(86,50)(86,140)(70,140)
\qbezier(85,50)(85,130)(70,130)
\qbezier(84,50)(84,120)(70,120)
\qbezier(83,50)(83,110)(70,110)
\qbezier(82,50)(82,100)(70,100)
\qbezier(81,50)(81,90)(70,90)
\qbezier(80,50)(80,80)(70,80)

\qbezier[55](70,40)(70,105)(70,150)

\put(47,30){$t=\pm 0.3$}

\put(160,50){\line(1,0){40}}
\put(160,60){\line(1,0){40}}
\put(160,70){\line(1,0){40}}
\put(160,80){\line(1,0){40}}
\put(160,90){\line(1,0){40}}
\put(160,100){\line(1,0){40}}
\put(160,110){\line(1,0){40}}
\put(160,120){\line(1,0){40}}
\put(160,130){\line(1,0){40}}
\put(160,140){\line(1,0){40}}

\put(300,50){\line(-1,0){40}}
\put(300,60){\line(-1,0){40}}
\put(300,70){\line(-1,0){40}}
\put(300,80){\line(-1,0){40}}
\put(300,90){\line(-1,0){40}}
\put(300,100){\line(-1,0){40}}
\put(300,110){\line(-1,0){40}}
\put(300,120){\line(-1,0){40}}
\put(300,130){\line(-1,0){40}}
\put(300,140){\line(-1,0){40}}

\qbezier(200,50)(214,50)(214,140)
\qbezier(200,60)(213,60)(213,140)
\qbezier(200,70)(212,70)(212,140)
\qbezier(200,80)(211,80)(211,140)
\qbezier(200,90)(210,90)(210,140)
\qbezier(200,100)(209,100)(209,140)
\qbezier(200,110)(208,110)(208,140)
 
\qbezier(260,50)(246,50)(246,140)
\qbezier(260,60)(247,60)(247,140)
\qbezier(260,70)(248,70)(248,140)
\qbezier(260,80)(249,80)(249,140)
\qbezier(260,90)(250,90)(250,140)
\qbezier(260,100)(251,100)(251,140)
\qbezier(260,110)(252,110)(252,140)

\qbezier(217,50)(217,140)(230,140)
\qbezier(218,50)(218,130)(230,130)
\qbezier(219,50)(219,120)(230,120)
\qbezier(220,50)(220,110)(230,110)
\qbezier(221,50)(221,100)(230,100)
\qbezier(222,50)(222,90)(230,90)
\qbezier(223,50)(223,80)(230,80)

\qbezier(243,50)(243,140)(230,140)
\qbezier(242,50)(242,130)(230,130)
\qbezier(241,50)(241,120)(230,120)
\qbezier(240,50)(240,110)(230,110)
\qbezier(239,50)(239,100)(230,100)
\qbezier(238,50)(238,90)(230,90)
\qbezier(237,50)(237,80)(230,80)

\linethickness{1pt}
\put(215.5,49){\line(0,1){92}}
\put(244.5,49){\line(0,1){92}}
\thinlines
 
\qbezier[55](230,40)(230,105)(230,150)

\qbezier(215.5,145)(215.5,160)(180,160)
\put(215.5,142){\vector(0,-1){0}}
\qbezier(244.5,145)(244.5,160)(180,160)
\put(244.5,142){\vector(0,-1){0}}
\put(132,158){{\small torus leaf}}

\put(216,30){$t=0$}

\end{picture}
\end{center}
\caption{Create torus leaf (2D movie)}
\label{fig-create.torus2}
\end{figure}
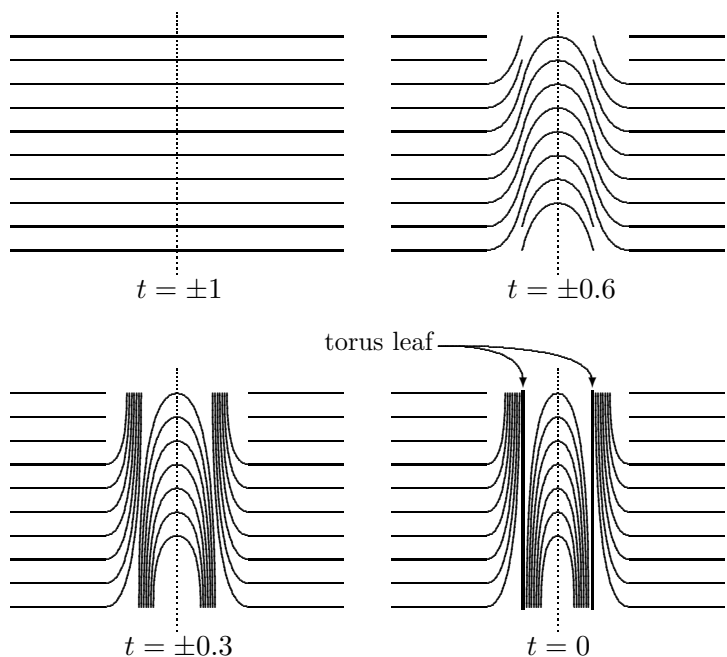
}

\section{Context}
\label{sec-context}

\subsection{Foliations and Gromov compactness}

First, an example that shows the flexibility of foliations:

\begin{example}[Y.\ Eliashberg]
\label{ex-create.torus}
\emph{Creating a torus leaf.}
Let  $\F$ be \emph{any} foliation on a $4$--manifold $M$. Let $c\co\Sph{1}\to M$ be any embedding. The curve $c$ can always be slightly perturbed to be transverse to $\F$. Choose another local coordinate near $c$, transverse both to $\F$ and to $c$, and think of it as \emph{time} (with $c$ appearing at time $t=0$). Start at time $t=-1$. As time goes on, begin pushing more and more the leaves of $\F$ parallel with the direction of $c$, wrapping them around more and more as time approaches $t=0$ 
(see Figures \ref{fig-create.torus1} and \ref{fig-create.torus2}).
At $t=0$, we can fit in a torus leaf, with the interior of the torus foliated by leaves diffeomorphic to $\aR^2$---a \defemph{Reeb component}. As time goes on from $t=0$, we play the movie backward. Notice that the new foliation is \emph{homotopic} with the one we started with(\ie the tangent plane fields are homotopic through integrable plane fields). (In particular, we have a geometric proof for part of \Corref{cor-torus}.)
\end{example}

\begin{corollary}
\label{cor-no.gromov}
Any \ac\ \str\ is homotopic to one for which Gromov compactness fails.
\end{corollary}

\emph{``Gromov compactness''} here means the compactness of the space of all $J$--holomorphic curves (curve = real surface). In other words, any sequence of holomorphic curves $f_n\co(\Sigma,j_n)\to (M,J)$ has a subsequence converging to a limit $f\co(\Sigma^*,j)\to(M,J)$ that is holomorphic (and may have nodal singularities; the limit domain $\Sigma^*$ is obtained by collapsing circles in $\Sigma$). Gromov compactness always holds for \ac\ \str s that admit symplectic \str s. (In fact, the essential property needed is that the areas of $f_n(\Sigma)$ be bounded; for a thorough discussion, see \cite{audin.lafontaine}.)

\begin{proof}
By \Thmref{thm-exist}, an \ac\ \str\ can be deformed till there is a singular foliation $\F$ with all leaves holomorphic. 
As in \Exref{ex-create.torus} above, create a torus leaf. 
Actually, by ``freezing'' the movie at $t=0$ (expanding the frame at $t=0$ to all $t\in[-\epsilon,\epsilon]$), create a lot of tori.
Now pick a second curve, orthogonal to these tori, and apply that example again.  What appears in the end is a torus that explodes to make room for a new Reeb component. Thinking in terms of an \ac\ \str\ $J$ induced by the final foliation, we have a sequence of $J$--holomorphic tori that has no decent limit.
\end{proof}

\begin{question}[R.~Kirby 2002]
\label{q-kirby}
Let $\F$ be a foliation on $M$, and $J$ an \ac\ \str\ making the leaves $J$--holomorphic. What conditions imposed on $\F$ insure that Gromov compactness holds for $J$?
\end{question}

In the extreme, if $\F$ is a Lefschetz pencil, then Gromov compactness holds (the manifold is symplectic). Compare also with \Propref{prop-2taut.sympl}.

The flexibility from our examples, at least, is done away with if we require the foliations to be \emph{taut} (since Reeb components kill tautness).

\subsection{Taut foliations}
\label{sec-taut}

A foliation $\F$ on a \riem\ manifold $(M,g)$ is called \defemph{minimal} if all its leaves are minimal surfaces in $(M,g)$ 
(\ie they locally minimize area; 
for any compact piece $K$ of a leaf, any small perturbation of $K$ rel $\del K$ will have bigger area; 
that is equivalent to each leaf having zero mean curvature).

A foliation $\F$ on $M$ is called \defemph{taut} if there is a \riem\ metric $g$ \st\ $\F$ is minimal in $(M,g)$.
(See  \cite[ch.~10]{candelconlon} for a general discussion.)

\begin{remark}
In the special case of a codimension-$1$ foliation, tautness is equivalent to the existence of a $1$--manifold transverse to $\F$ and crossing all the leaves. 
A similar condition is too strong for higher codimensions.
\end{remark}
 
Tautness can be expressed in terms of $2$--forms \cite{rummler}:

\begin{theorem}[H.~Rummler 1979]
\label{thm-rummler}
Let $\F$ be a foliation on $M$. Then $\F$ is taut \Iff\ there is a $2$--form $\mu$ \st\ $\mu\rest{\text{Leaf}}>0$ and
$d\mu\rest{\text{Leaf}}=0$.
\end{theorem}

We write 
``$\mu\rest{\text{Leaf}}>0$ 
and $d\mu\rest{\text{Leaf}}=0$'' 
as shorthand for 
``$\mu(\tau_1,\tau_2)>0$ and $d\mu(\tau_1,\tau_2,z)=0$, for any orienting pair $\tau_1, \tau_2\in\T{\F}$ and any $z\in\T{M}$''.

The strong link between $2$--forms and minimality of foliations is also suggested (via \ac\ \str s) by the following formula:

\begin{lemma}
\label{lemma-domega}
Let $g(x,y)=\inner{x,y}$ be a \riem\ metric on $M^{4}$ and $\nabla$ 
its \LC\ connection. Let $J$ be any $g$--orthogonal \ac\ \str, and
let $\omega(x,y)=\inner{Jx,\,y}$ be its fundamental $2$--form.
Let $x$, $z$ be any vector fields on $M$. Then:
\[ d\omega(x,Jx,z)
 =\Inner{[x,Jx],\,Jz}
 -\Inner{\Nabla{x}x+\Nabla{Jx}Jx,\ z} \]
\end{lemma}

The term $[x,Jx]$ measures the integrability of the $J$--holomorphic plane field $\aR\inner{x,Jx}$, while the normal component of the term $\Nabla{x}x+\Nabla{Jx}Jx$ is the mean curvature of the plane field $\aR\inner{x,Jx}$, and thus measures its $g$--minimality. (\Lemmaref{lemma-domega} will be proved at the end of the paper, in \Secref{sec-appendix}.)

\subsection{Minimal genus of embedded surfaces}
\label{sec-min.genus}

Given any class $a\in H_2(M;\Zi)$, there are always embedded surfaces in $M$ that represent it. An open problem is to determine how simple such surfaces can be, or, in other words, what is the minimal genus that a surface representing $a$ can have. (Remember that $\chi(S)=2-2g(S)$, so minimum genus is maximum Euler characteristic.)

Notice that, if $S$ is a $J$--holomorphic surface for some \ac\ \str\ $J$, then 
\[ \chi(S)+S\cdot S=c_1(J)\cdot S \]
(simply because $c_1(J)\cdot S=c_1(\T{M}\rest{S})=c_1(\T{S})+c_1(\N{S})=\chi(S)+S\cdot S$). 
This equality is known as the ``adjunction formula'' for $S$.

In general, the main and most powerful tool for obtaining genus bounds comes from \SW\ theory \cite{kronheimer.mrowka.proj,ozsvath.szabo.sympl.thom}:

\begin{proposition}[Adjunction Inequality]
Let $S$ be any embedded surface in $M$. Assume that either $M$ is of \SW\ simple type and $S$ has no sphere components, or that $S\cdot S\geq 0$.
Then, for any \SW\ basic class $\boldsymbol{\epsilon}$, we have:
\[ \chi(S)+S\cdot S \leq 
	\boldsymbol{\epsilon}\cdot S \]
In particular, if $J$ is an \ac\ \str\ admitting a symplectic structure, then 
\[ \chi(S)+S\cdot S \leq c_1(J)\cdot S \] 
\end{proposition}

Nonetheless, the bounds offered by \SW\ basic classes are not always sharp.  
For example, in the case of manifolds $M=N^3\times\Sph{1}$, P.~Kronheimer has proved in \cite{kronheimer} that foliations give \emph{better} bounds:

\begin{theorem}[P.\ Kronheimer 1999]
\label{thm-kron}
Consider $M=N^3\times\Sph{1}$, with $N$ a closed irreducible $3$--manifold.
Let $\Bar{\F}$ be a taut foliation in $N$, with Euler class $\bar{\boldsymbol{\epsilon}}=e(\T{\Bar{\F}})$. Let $\boldsymbol{\epsilon}$ be the image of $\bar{\boldsymbol{\epsilon}}$ in $H^2(N\times\Sph{1})$. Then, for any embedded surface $S$ in $M$, without sphere components, we have
\[ \chi(S)+S\cdot S \leq 
	\boldsymbol{\epsilon}\cdot S \]
\end{theorem}
In general, $\boldsymbol{\epsilon}$ is \emph{not} a \SW\ basic class. (Nonetheless, the proof of \Thmref{thm-kron} does use \SW\ theory: the taut foliation $\Bar{\F}$ is perturbed to a tight contact structure, which is then symplectically filled in a suitable way, and a version of the \SW\ invariants is used: $\boldsymbol{\epsilon}$ is a ``monopole class''.)

Taut foliations $\Bar{\F}$ on $3$--manifolds are well-understood, and are stron\-gly related to minimal genus surfaces there:
If $N^3$ is a closed irreducible $3$--manifold, then an embedded surface $S$ has minimal genus \Iff\ it is the leaf of a taut foliation of $N$
\cite{thurston-norm,gabai}. 
(A similar statement on $4$--manifolds is not known.)

A taut foliation $\Bar{\F}$ on $N^3$ induces an obvious product foliation $\F=\Bar{\F}\times\Sph{1}$ on $M=N\times\Sph{1}$ (with leaves $\text{Leaf}\times\{pt\}$). Then $\F$ is also  taut (pick a product metric), $\boldsymbol{\epsilon}=e(\T{\F})$ is the pull-back of $\bar{\boldsymbol{\epsilon}}=e(\T{\Bar{\F}})$, and
the \ac\ \str\ that $\F$ determines has $c_1(J_\F)=\boldsymbol{\epsilon}=e(\T{\F})$. One can then try to generalize \Thmref{thm-kron} as

\begin{conjecture}[P.\ Kronheimer 1999]
\label{conj}
Let $\F$ be a taut foliation on $M^4$, and $J_\F$ be an \ac\ \str\ induced by $\F$. Then, for any embedded surface $S$ without sphere components, we have
\[ \chi(S)+S\cdot S\leq c_1(J_\F)\cdot S \]
\end{conjecture}

\noindent
A few extra requirements are needed, \eg to exclude manifolds like $S\times\Sph{2}$. Kronheimer also proposes that the foliations be allowed singularities.

\begin{remark}
The situation in \Thmref{thm-kron} has another peculiarity: $\F$ admits  transverse foliations. Indeed, since $\Bar{\F}$ has codimension $1$, any nowhere-zero vector-field in $N^3$ normal to $\Bar{\F}$ integrates to a $1$--dimensional foliation of $N$ that is transverse to $\Bar{\F}$. By multiplying its leaves by $\Sph{1}$, this $1$--dimensional foliation induces a $2$--dimensional foliation in $M$ that is transverse to $\F$.
\end{remark}

One could thus think of strengthening the hypothesis of \Conjref{conj} by requiring $\F$ not only to be taut, but also to admit a transverse foliation.
One might push things even further and ask that the second foliation be taut as well. But then one almost runs into:

\begin{proposition}[Two taut makes one symplectic]
\label{prop-2taut.sympl}
Let $\F$ and $\G$ be transverse foliations on $M^4$. If there is a metric $g$ that makes both $\F$ and $\G$ be minimal and orthogonal, then $M$ must admit a symplectic \str. 
Therefore, for any embedded surface $S$ we have
\[ \chi(S)+S\cdot S\leq c_1(J_\F)\cdot S \]
\end{proposition}

\noindent
(This is an immediate consequence of \Lemmaref{lemma-domega}.)

\begin{remark}
\emph{No taut, no symplectic.}
At the other extreme, if $M$ admits a \emph{non-taut} foliation $\F$, then \emph{no} \ac\ \str\ compatible with $\F$ admits symplectic \str s (or, more directly, $M$ admits no symplectic structures making the leaves of $\F$ symplectic submanifolds).
\end{remark}


\section{Proofs}
\label{sec-proofs}

\subsection{Thurston's theorem}

The tool that we use for obtaining foliations is the $h$--principle discovered by W.~Thurston \cite{thurston} for such objects:

\begin{theorem}[W.\ Thurston 1974]
\label{thm-thurston}
Let ${\cal T}$ be a $2$--plane field on a manifold $M$ of dimension at least $4$. Let $K$ be a compact subset of $M$ \st\ ${\cal T}$ is completely integrable in a neighborhood of $K$ ($K$ can be empty). Then ${\cal T}$ is homotopic $\operatorname{rel} K$ to a completely integrable plane field.
\end{theorem} 

This theorem is also true on $3$--manifolds (see \cite{thurston-codim1}), but \emph{not} in a relative version.
The theorem above is proved by first lifting the plane field to a Haefliger structure, and then deforming the latter to become a foliation using the main theorem of \cite{thurston}. The latter result has an alternative proof in \cite{eliashberg}.

\begin{remark}
A problem with using Thurston's theorem is that, when following its proof to build foliations, one only gets \emph{non-taut} foliations. Indeed, certain holes in the foliation being built have to be filled-in with Reeb components.
\end{remark}

Thurston's theorem reduces the problem of building foliations to the problem of finding singular plane-fields on $M$, or, more exactly, singular splittings $\T{M}=T\oplus N$. ``Singular'' because the difference between $\T{M}$ and the sum $T\oplus N$ will be a surgery modification that we present next:

\subsection{Surgery modifications}

If $E\to M$ is an oriented $4$--plane bundle and $B$ is a $4$--ball around a point $x$, then we can cut out $E\rest{B}$ and glue it back in using an automorphism of 
$E\rest{\del B}$.
Choose a chart in $M$ around $B$ and use some quaternion coordinates $\aR^{4}\iso\Ha$, so that $\del B\iso\Sph{3}$, the sphere of units in $\Ha$. 
Since the fiber of $E$ is $4$--dimensional and $E\rest B$ is trivial, we can choose some quaternion bundle-coordinates on $E$, so that $E\rest{B}\iso\Ha\times B$. Then $E\rest{\del B}\iso\Ha\times\Sph{3}$.

Quaternions can be used to represent $SO(4)$ acting on $\aR^4$ as $\Sph{3}\times\Sph{3}\big/\pm1$ acting on $\Ha$ through $[q_+,q_-]h=q_+ h q_-^{-1}$.
For any $m,n\in\Zi$, we define a map $\Sph{3}\to SO(4)$ by 
\[ \xi_{m,n}(q)h=q^m	h q^n \]	
where $q\in\Sph{3}$ and $h\in\Ha\iso\aR^4$.
Notice that the map $\xi_{m,n}$ determines an element of $\pi_3 SO(4)$. In fact, we have the isomorphism $\Zi\oplus\Zi\iso\pi_3 SO(4)$ given by $(m,n)\to[\xi_{m,n}]$ (see \cite{steenrod}). Homotopically we have $[\xi_{m,n}]+[\xi_{p,q}]=\xi_{m,n}\circ\xi_{p,q}=[\xi_{m+p,\,n+q}]$.

The $(m,n)$--surgery modification of $E$ is then defined as follows: Pick a point $x$ in $M$ and a $4$--ball $B$ around it. Cut $E\rest{B}$ out from $E$ and glue it back using the automorphism 
\[\begin{CD}
	E\rest{\del B} @>>> E\rest{\del B}\qquad\\
	(x,h) &\maps& \bigl(x,\ \xi_{m,n}(x)h\bigr)
\end{CD}\]
Denote the resulting bundle by $E_{m,n}$.

\begin{remark}
It is equivalent to perform $(m,n)$--surgery at one point $x$, or to perform $(1,0)$--surgeries at $m$ points $x_1,\ldots,x_m$ and $(0,1)$--surgeries at $n$ points $x'_1,\ldots,x'_n$ (or any other combination that adds up to $(m,n)$ in $\Zi\oplus\Zi$).
\end{remark}

To understand the result of such a modification, we study its characteristic classes:

\subsection{Characteristic classes}

An oriented $k$--bundle over $\Sph{n}$ is unique\-ly determined by the homotopy class of an equatorial gluing map $\Sph{n-1}\to SO(k)$, and thus $\operatorname{Vect}_k \Sph{n}\iso\pi_{n-1} SO(k)$.

In particular, all oriented $4$--bundles on $\Sph{4}$ correspond one-to-one with $\pi_3 SO(4)$. Therefore all of them can be obtained by $(m,n)$--surgery modifications. Denote by $\uaR^4_{m,n}$ the bundle on $\Sph{4}$ obtained from the trivial bundle $\uaR^4=\aR^4\times\Sph{4}$ through a $(m,n)$--modification.

Note that addition of gluing maps in $\pi_3 SO(4)$ survives as addition of characteristic classes of bundles in $H^*(\Sph{4};\Zi)$.
In particular, for any characteristic class $\goth{c}$, we have 
\[ \goth{c}(\uaR^4_{m,n})
	=m\goth{c}(\uaR^4_{1,0})+n\goth{c}(\uaR^4_{0,1}) \] 

It is known that $(1,1)$--surgery on $\aR^4\times\Sph{4}$ will yield the tangent bundle $\T{\Sph{4}}$ of $\Sph{4}$ (see \cite{steenrod}). Since $e(\T{\Sph{4}})=\chi(\Sph{4})=2$, we deduce that 
\[ e(\uaR^4_{1,0})+e(\uaR^4_{0,1})=2 \]
On the other hand, $\T{\Sph{4}}\oplus\uaR=\uaR^5$, so $p_1(\T{\Sph{4}})=p_1(\T{\Sph{4}}\oplus\uaR)=0$, and so
\[ p_1(\uaR^4_{1,0})+p_1(\uaR^4_{0,1})=0 \]

The bundle $\uaR^4_{1,-1}$ is obtained by surgery with $\xi_{1,-1}(q)h=qhq^{-1}$. The latter preserves the real line $\aR\subset\Ha$ of the fiber, and therefore the bundle $\uaR^4_{1,-1}$ splits off a trivial real-line bundle. Therefore $e(\uaR^4_{1,-1})=0$. That means 
\[ e(\uaR^4_{1,0})-e(\uaR^4_{0,1})=0 \]
Combining with the above yields
\[ e(\uaR^4_{1,0})=1  \qquad\qquad e(\uaR^4_{0,1})=1 \]

\begin{remark}
\label{rk-cx.quaternion}
\emph{Complex structures on quaternions}\qua
The complex plane $\Ce^2$ can be identified with the quaternions $\Ha$ in two ways: 
\begin{enumerate} 
\item 
$(z_1,z_2)\equiv z_1+z_2 j$ (and then complex scalars are multiplying in $\Ha$ on the left, and the natural orientations of $\Ce^2$ and $\Ha$ are preserved; quaternion multiplication on the right is $\Ce$--linear, and $\Sph{3}$ acting on the right identifies with $SU(2)$)
\item
$(z_1,z_2)\equiv z_1+jz_2$ (with complex scalars multiplying on the right, but with the orientations \emph{reversed}; quaternions multiplying on the left act $\Ce$--linearly, $\Sph{3}$ on the left is $SU(2)$). 
\end{enumerate}
(We will make use of both of these identifications: (1) will be used here, while (2) in \S\ref{sec-not.cx}.)
\end{remark}

We identify $\Ce^2$ with $\Ha$ using $(z_1,z_2)\equiv z_1+z_2 j$. Since the bundle $\uaR^4_{0,1}$ is built using the map $\xi_{0,1}(q)h=hq$, and the latter preserves multiplication by complex scalars on the left, we deduce that $\uaR^4_{0,1}$ can be seen as a complex-plane bundle (the same is true for all $\uaR^4_{0,n}$).
Thus $\uaR^4_{0,1}$ has well-defined Chern classes. Since $c_2(\uaR^4_{0,1})=e(\uaR^4_{0,1})$ and $c_1(\uaR^4_{0,1})\in H^2(\Sph{4})$, we see that
\[ c_1(\uaR^4_{0,1})=0 \qquad\qquad c_2(\uaR^4_{0,1})=1 \]
Since for complex bundles we have $p_1=c_1^2-2c_2$, we deduce that
\[ p_1(\uaR^4_{0,1})=-2 \]
and therefore, combining with the above, 
\[ p_1(\uaR^4_{1,0})=2 \]
Therefore:
\[ e(\uaR^4_{m,n})=m+n \qquad\qquad p_1(\uaR^4_{m,n})=2m-2n \]

In conclusion, 
for any oriented $4$--plane bundle $E\to\Sph{4}$, we have 
\[ e(E_{m,n})
	=e(E)+m+n \qquad\qquad p_1(E_{m,n})=p_1(E)+2m-2n \]
This change of characteristic classes for bundles over $\Sph{4}$ is also what happens over a general $4$--manifold $M$.  This can be seen, for example, using the obstruction-theoretic definition of characteristic classes (defined locally cell-by-cell): away from the modification, it does not matter if we are left with a small neighborhood of the south pole, or with $M\setminus\text{\it Ball}$.
(Or, one could argue that $\bigl[M\setminus\text{\it Ball},\ BSO(4)\bigr]$ is finite, while $p_1$ are $e$ are rational, etc.)

\begin{lemma}
For any $4$--plane bundle $E\to M$, a $(m,n)$--modification of $E$ will change its characteristic classes as follows:
\[ e(E_{m,n})=e(E)+m+n \qquad\qquad p_1(E_{m,n})=p_1(E)+2m-2n \]
\end{lemma}

\subsection{Obtaining the tangent bundle}

Assume $c\in H^2(M;\Zi)$ is an integral lift of the Stiefel--Whitney class $w_2(M)\in H^2(M;\Zi_2)$. 
For any splitting $c=\tau+\nu$, build the complex-line bundles $L_\tau$ and $L_\nu$ such that $c_1(L_\tau)=\tau$ and $c_1(L_\nu)=\nu$. Let $E=L_\tau\oplus L_\nu$. Then $c_1(E)=c$ and $c_2(E)=\tau\nu$. Thus, as a real $4$--plane bundle, $E$ has $w_2(E)=w_2(M)$. 

If we can modify $E$ to an $E'$ so that we also have $e(E')=\chi(M)$ and $p_1(E')=p_1(M)$, then $E'\iso\T{M}$. 
That is due the fact that characteristic classes determine bundles up to isomorphism \cite{doldwhitney}:

\begin{theorem}[A.\ Dold \& H.\ Whitney 1959]
\label{thm-dold.whitney}
Let $E_1\to M$ and $E_2\to M$ be two oriented $4$--plane bundles over an oriented $4$--manifold $M$. If $w_2(E_1)=w_2(E_2)$, $e(E_1)=e(E_2)$, and $p_1(E_1)=p_1(E_2)$, then $E_1\iso E_2$.
\end{theorem}

Now, since 
\begin{align*}
e(E_{m,n}) &=e(E)+m+n=\tau\nu+m+n \\ 
p_1(E_{m,n}) &=p_1(E)+2m-2n=c_1(E)^2-2c_2(E)+2m-2n\\
  &=c^2-2\tau\nu+2m-2n
\end{align*}
we obtain that $e(E_{m,n})=e(\T{M})$ and $p_1(E_{m,n})=p_1(\T{M})$
\Iff\
\begin{align*}
& m=\rec{4}\bigl( p_1(M)+2\chi(M)-c^2 \bigr)\\
& n=\rec{4}\bigl( -p_1(M)+2\chi(M)+c^2-4\tau\nu \bigr)
\end{align*}

\begin{remark}
These $m$ and $n$ are \emph{always integers}. The quick argument is:
on the one hand, the formula for $m$ above gives exactly the dimension of the \SW\ moduli space associated to the \spinc-\str\ given by $c$ (it is the index of a differential operator), and thus is known to be integral; on the other hand, $n=-m+\chi(M)-\tau\nu$.
\end{remark}

In conclusion:

\begin{proposition}[Splitting the Tangent Bundle]
\label{prop-kirby}
Let $\tau, \nu\in H^2(M;\Zi)$ be such that $c=\tau+\nu$ is an integral lift of $w_2(M)\in H^2(M;\Zi_2)$.
Let $L_\tau, L_\nu$ be complex-line bundles with $c_1(L_\tau)=\tau$ and $c_1(L_\nu)=\nu$.
Then
\[ (L_\tau\oplus L_\nu)_{m,n}\iso\T{M} \]
where
\begin{align*}
& m=\rec{4}\bigl( p_1(M)+2\chi(M)-c^2 \bigr)\\
& n=\rec{4}\bigl( -p_1(M)+2\chi(M)+c^2-4\tau\nu \bigr)
\end{align*}
\end{proposition}

In the case $\nu=0$, this is a statement that we learned (together with its proof) from R.~Kirby's lectures at U.~C.~Berkeley. (The advantage of using a more complicated sum $L_\tau\oplus L_\nu$ versus the simpler $L_c\oplus\uaR^2$ will become apparent when we move toward foliations.)

In the special case when $c$ is a complex class 
(\ie when $c$, besides being an integral lift of $w_2(M)$, also satisfies
$p_1(M)=c^2-2\chi(M)$),
we have 
\[ m=0 \qquad\qquad n=\chi(M)-ab \]
Since $m=0$, that means, in particular, that the surgery modification is made with $\xi_{0,n}(q)h=hq^n$, which is $\Ce$--linear, and thus will preserve the complex structure of $L_\tau\oplus L_\nu$. Therefore $\T{M}$ inherits a complex-structure. We have thus built an \ac\ \str\ $J$ on $M$ with $c_1(\T{M},J)=c$.

\begin{corollary}
\label{cor-cx.str}
Let $\tau, \nu\in H^2(M;\Zi)$ be such that $\tau+\nu$ is a complex class.
Let $L_\tau, L_\nu$ be complex-line bundles with $c_1(L_\tau)=\tau$ and $c_1(L_\nu)=\nu$.
Then $M$ admits an \ac\ \str\ $J$ with $c_1(J)=\tau+\nu$, and, for $n=\chi(M)-\tau\nu$, we have
\[ (\T{M},J)\iso(L_\tau\oplus L_\nu)_{0,n} \]
as complex bundles.
\end{corollary}

\begin{remark}
If $\chi(M)-\tau\nu\geq 0$, then the complex bundle $(\T{M},J)$ can be obtained from $L_\tau\oplus L_\nu$ by surgery modifications at $n=\chi(M)-\tau\nu$ points using $\xi_{0,1}$.
\end{remark}

Through the isomorphism $(L_\tau\oplus L_\nu)_{m,n}\iso\T{M}$, the line-bundles $L_\tau$ and $L_\nu$ survive as plane-fields $\til{L}_\tau$ and $\til{L}_\nu$ in $\T{M}$ defined off the modification points. If we find a way to prolong $\til{L}_\tau$ across its singularities by a singular foliation, then we could use Thurston's \Thmref{thm-thurston} (in its relative version) to integrate the whole $\til{L}_\tau$ to a foliation $\F$ (while keeping it fixed at the singularities). This foliation would then, off the singularities, have $\T{\F}\iso L_\tau$ and $\N{\F}\iso L_\nu$, and thus have well-defined Euler classes $e(\T{\F})=\tau$ and $e(\N{\F})=\nu$ in $H^2(M;\Zi)$ (since the isolated singular points cannot influence $H^2$, see \Rkref{rk-euler}). 

Finding nice singularities is what we do next:

\subsection{Singularities}

We keep identifying $\Ha\iso\Ce^2$ using $z_1+z_2 j\equiv(z_1,z_2)$. In particular $\aR\inner{1,i}=\Ce$ in $\Ha$ is $\Ce\times{0}$ in $\Ce^2$.

Consider the action of $\xi_{0,1}$ on $\Sph{3}\times\Ha$. Since $\xi_{0,1}(q)\cdot 1=q$ and $\xi_{0,1}$ is $\Ce$--linear, we deduce that $\xi_{0,1}(q)\cdot \Ce =\Ce q$. 
In other words, the trivial subbundle $\Sph{3}\times\Ce$ is taken by $\xi_{0,1}$ to the subbundle $\bigcup\, \{q\}\times\Ce q$ whose fiber over a point $q$ of $\Sph{3}$ is the complex plane spanned by $q$.

Assume now that $c=\tau+\nu$ is a complex class and that $n=\chi(M)-\tau\nu \geq 0$. Then we can  build $\T{M}$ from $L_\tau\oplus L_\nu$ as above, by modifying at $n$ points using $\xi_{0,1}$. 

Choose coordinates on a small $4$--ball $B$ around a modification point $x$ so that the fibers of $L_\tau$ on $\Sph{3}=\del B$ are $\Ce\subset\Ha$. Then $\xi_{0,1}$ will glue the fiber of $L_\tau$ over $q\in\Sph{3}$ to the plane $\Ce q$. The latter can be identified though with the tangent planes to the submanifolds $\Ce q$ of the unit ball in $\Ce^2$ bounded by $\Sph{3}$. Or, in other words, $q\maps \xi_{0,1}(q)\cdot L_\tau$ is tangent to the levels of the complex function $(z_1,z_2)\maps z_1/z_2$ (from $\Ce^2\setminus 0$ to $\CP^1$).

\begin{figure}[t] 
\begin{center}
\setlength{\unitlength}{0.9pt}
\begin{picture}(360,220)(-10,-60)

\put(50,100){\bigcircle{100}}

\put(-10,100){\line(1,0){20}} 
\put(35,150){\line(1,0){30}}
\put(90,100){\line(1,0){20}}
\put(35,50){\line(1,0){30}}

\put(2,66){\line(1,0){22}} 
\put(98,66){\line(-1,0){22}}
\put(98,134){\line(-1,0){22}}
\put(2,134){\line(1,0){22}}

\put(-6,83){\line(1,0){21}} 
\put(-6,117){\line(1,0){21}}
\put(106,83){\line(-1,0){21}}
\put(106,117){\line(-1,0){21}}

\put(13,143){\line(1,0){23}}
\put(13,56){\line(1,0){23}}
\put(87,143){\line(-1,0){23}}
\put(87,56){\line(-1,0){23}}

\put(-10,40){\line(1,0){20}} 
\put(-10,150){\line(1,0){20}}
\put(110,40){\line(-1,0){20}}
\put(10,15){\line(1,0){20}}
\put(60,-20){\line(1,0){20}}
\put(-10,-60){\line(1,0){20}}

\put(210,100){\bigcircle{100}}

\put(150,100){\line(1,0){20}} 
\put(210,160){\line(0,-1){20}}
\put(270,100){\line(-1,0){20}}
\put(210,40){\line(0,1){20}}

\put(167.5,57.5){\line(1,1){14}} 
\put(167.5,142.5){\line(1,-1){14}}
\put(252.5,142.5){\line(-1,-1){14}}
\put(252.5,57.5){\line(-1,1){14}}

\put(155,122){\line(3,-1){18}} 
\put(155,78){\line(3,1){18}}
\put(265,122){\line(-3,-1){18}}
\put(263,78){\line(-3,1){18}}

\put(190,157){\line(1,-3){6.5}}
\put(230,157){\line(-1,-3){6.5}}
\put(190,43){\line(1,3){6.5}}
\put(230,43){\line(-1,3){6.5}}

\put(290,0){\bigcircle{100}}

\put(230,0){\line(1,0){120}}
\put(290,-60){\line(0,1){120}}
\put(247.5,-42.5){\line(1,1){85}}
\put(247.5,42.5){\line(1,-1){85}}

\put(290,0){\line(3,1){56}} 
\put(290,0){\line(-3,1){56}}
\put(290,0){\line(3,-1){56}}
\put(290,0){\line(-3,-1){56}}

\put(290,0){\line(1,3){19}}
\put(290,0){\line(-1,3){19}}
\put(290,0){\line(1,-3){19}}
\put(290,0){\line(-1,-3){19}}

\thicklines
\qbezier(110,125)(130,130)(150,125)
\put(150,125){\vector(4,-1){0}}
\put(125,135){\large $\xi_{0,1}$}

\qbezier[50](177,57)(225,28)(243,-32)
\qbezier[50](253,133)(310,100)(335,35)

\put(185,-40){$z_1/z_2=\epsilon$}

\end{picture}
\end{center}
\caption{Filling with a pencil singularity.}
\label{fig-pencil}
\end{figure}
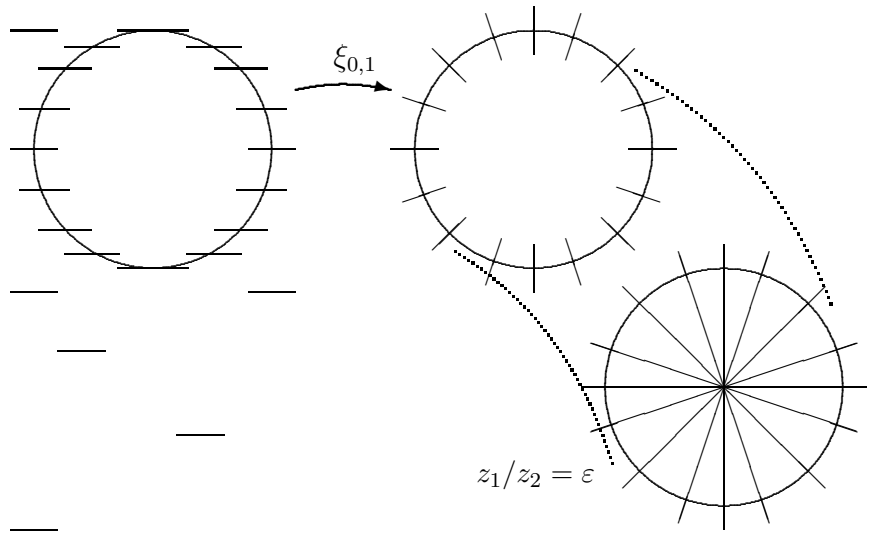

These levels can be used to fill-in the singularity of the foliation $\F$ obtained by deforming $L_\tau$ (see \Figref{fig-pencil}). We call such a (filled-in) singularity a \defemph{pencil singularity}.

\medskip

Since we are dealing with bundles, though, what essentially matters when filling a singularity is only the \emph{homotopy class} of its boundary plane field $q\maps \xi_{0,1}(q)\cdot L_\tau$, seen as a map $\Sph{3}\to\CP^1$. 
That is completely determined by the homotopy class 
of any spanning vector field
(for example $q\maps\xi_{0,1}(q)\cdot 1=q$),
seen as a map $\Sph{3}\to\Sph{3}$. 

\begin{remark}
The two homotopy classes are related by the Hopf map $\goth{h}:\Sph{3}\to\CP^1=\Sph{2}$, which establishes the isomorphism $\pi_3\Sph{3}\iso\pi_3\Sph{2}$. Technically, a map $u\co\Sph{3}\to\Sph{3}$ has a \emph{degree}, while a map $v\co\Sph{3}\to\CP^1$ has a 
\emph{Hopf invariant}. 
When $v=\Ce\cdot u$ (that is: $v=\goth{h}u$), the two coincide, and we will call them ``degree'' in both instances.
\end{remark}

Consider the levels of the complex function $(z_1,z_2)\maps z_1 z_2$. The tangent space to the level through $(z_1,z_2)$ is the complex span of $(z_1,-z_2)$. The latter, restricted to a map $\Sph{3}\to\Sph{3}$, has degree $1$. The plane field $\til{L}_a\rest{\Sph{3}}$ also has degree 1 (since it is spanned by $q\mapsto q$).
Therefore $\til{L}_a$ can be homotoped to become tangent to the levels of the function $(z_1,z_2)\maps z_1 z_2$.

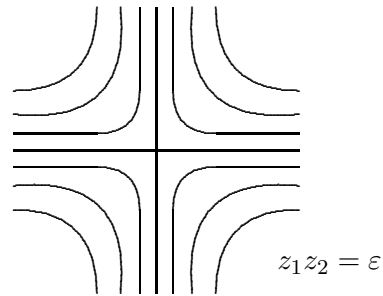
\begin{figure}[t] 
\begin{center}
\begin{picture}(120,120)(-10,0)


\thicklines
\put(-10,50){\line(1,0){120}}
\put(50,110){\line(0,-1){120}}
\thinlines

\put(-10,57){\line(1,0){35}}
\qbezier(25,57)(43,57)(43,75)
\put(43,75){\line(0,1){35}}

\put(110,57){\line(-1,0){35}}
\qbezier(75,57)(57,57)(57,75)
\put(57,75){\line(0,1){35}}

\put(-10,43){\line(1,0){35}}
\qbezier(25,43)(43,43)(43,25)
\put(43,25){\line(0,-1){35}}

\put(110,43){\line(-1,0){35}}
\qbezier(75,43)(57,43)(57,25)
\put(57,25){\line(0,-1){35}}

\qbezier(25,110)(25,75)(-10,75)
\qbezier(75,110)(75,75)(110,75)
\qbezier(25,-10)(25,25)(-10,25)
\qbezier(75,-10)(75,25)(110,25)

\qbezier(35,110)(40,60)(-10,65)
\qbezier(65,110)(60,60)(110,65)
\qbezier(35,-10)(40,40)(-10,35)
\qbezier(65,-10)(60,40)(110,35)

\put(100,0){$z_1 z_2=\epsilon$}

\end{picture}
\end{center}
\caption{A quadratic singularity (fake image).}
\label{fig-quadratic}
\end{figure}

Thus the levels of $z_1 z_2$ offer another possible way of filling-in the singularities of $\F$ (see \Figref{fig-quadratic}). (Notice that these levels are isomorphic to the levels of $(z_1,z_2)\maps z_1^2+z_2^2$.) We will call such a (filled-in) singularity a \defemph{quadratic singularity}.

\begin{Proof}{Existence \Thmref{thm-exist}}
Let $c$ be a complex class of $M$, and $c=\tau+\nu$ a splitting such that $n=\chi(M)-\tau\nu$ is non-negative. Then, by \ref{cor-cx.str}, we have $(L_\tau\oplus L_\nu)_{0,n}\iso\T{M}$, for complex-line bundles $L_\tau$ and $L_\nu$ with $c_1(L_\tau)=\tau$ and $c_1(L_\nu)=\nu$. We choose to perform the surgery by  $n$ modifications by $\xi_{0,1}$ at $n$ points $p_1,\ldots,p_n$. The bundle $L_\tau$ survives in $\T{M}$ as a singular plane-field $\til{L}_\tau$. Choose any assortment of $n$ pencil or quadratic singularities, and place them at the points $p_1,\ldots,p_n$. Arrange $\til{L}_\tau$ so that it is tangent to the leaves of the singularities. Use Thurston's \thmref{thm-thurston} to homotop $\til{L}_\tau$ (away from the singularities) so that it becomes integrable. The resulting singular foliation $\F$ is what we needed to build.
\end{Proof}

Of course, many other singularities can be chosen, see \Secref{sec-other.sing} below.
We singled out $z_1/z_2$ and $z_1 z_2$ now because they are exactly the singularities that appear in a Lefschetz pencil. 
Unlike a Lefschetz pencil, though, \emph{not} all leaves must pass through a pencil singularity. (For example, use the method of \Exref{ex-create.torus} to create a torus leaf that does not touch any singularity).








The choice being given, pencil singularities are the most manageable:  A pencil singularity  can be removed by blowing $M$ up: the blow-up simply separates the leaves of $\F$ that were meeting there, and thus the foliation survives with one less singularity. That is not the case for a quadratic singularity: blowing-up creates one more singularity for the foliation (since the exceptional sphere must now become a leaf), and instead of one there are now two quadratic singularities. On the other hand, a quadratic singularity creates at most two singular leaves, while a pencil singularity creates uncountably many. Also, in rare occasions, if one of the leaves that passes through a quadratic singularity is a sphere with self-intersection $-1$, one might attempt to blow it down while preserving the rest of the foliation. (Notice that the existence of a sphere leaf with non-zero self-intersections is not excluded by Reeb stability if the leaf passes through a singularity.)

\subsection{Other singularities}
\label{sec-other.sing}

Other singularities may be used, as stated
in \Thmref{thm-exist.sing}.
What matters is their \emph{Hopf degree}, defined as follows: 

If $\F$ is a foliation with an isolated singularity at $p$, then choose a small $4$--ball $B$ around $p$. If the plane-field $\T{\F}\rest{\del B}$ is left invariant by some local \ac\ \str\ on $B$ that is compatible with the  orientation of $M$, then we call $p$ a singularity of \defemph{positive type}. 

For a singularity of positive type, we define its \defemph{Hopf degree} as the degree (Hopf invariant) of the plane field $\T{\F}\rest{\del B}$ seen as maps $\del B\to\CP^1$ (\ie $\Sph{3}\to\Sph{2}$).
(Technically, if one wants the Hopf degree to depend only on the singularity and not on the chosen neighborhood $B$, one should define the Hopf degree as a \emph{limit} as $B$ shrinks to $p$.)

The most obvious candidates for singularities of positive type are, of course, singularities defined by levels of complex functions. But since a complex function will always preserve orientations, \emph{any} singularity coming from a complex function will have \emph{non-negative} Hopf degree. Thus, while we can easily find singularities for any positive Hopf degree (that can be used to fill-in the singularities created by $(0,n)$--modifications when $n\geq 0$), singularities of negative Hopf degree seem harder to find. Thus, the author does not know how to fill-in the singularity created by $\xi_{0,-1}$, which is why many statements have a positivity condition like $\chi(M)-\tau\nu\geq 0$.

\begin{example}
\emph{Cusp singularity.}
The singularity defined by the levels $\F(f)$ of the function
\[ f(z_1,z_2)=z_1^3-z_2^2 \]
It has $Df\rest{z}=\bigl[3z_1^2,\, -2z_2\bigr]$, and thus $\T{\F(f)}\rest{z}=\ker Df\rest{z}=\Ce(2z_2,\, 3z_1^2)$. A generic equation $(2z_2,\,3z_1^2)=(w^0_1,w^0_2)$ has two solutions, orientations are preserved, and thus
\[ \deg\F(f)=2 \]
\end{example}

\begin{example}
\emph{Normal crossing.}
The singularity defined by 
\[ g(z_1,z_2)=z_1^p z_2^q \]
has $Dg\rest{z}=\bigl[pz_1^{p-1}z_2^q,\ qz_1^p z_2^{q-1}\bigr]$, so $\T{\F(g)}\rest{z}=\Ce(qz_1, -pz_2)$, and thus
\[ \deg\F(g)=1 \]
\end{example}


\begin{example}
The singularity defined by
\[ k(z_1,z_2)=z_1^{p+1}+z_2^{q+1}\]
has $Dk\rest{z}=\bigl[(p+1)z_1^p,\ (q+1)z_2^q\bigr]$, so
$\T{\F(k)}\rest{z}=\Ce\bigl((q+1)z_2^q,\, (p+1)z_1^p\bigr)$, and thus
\[ \deg\F(k)=pq \]
\end{example}

This last example shows that all positive degrees are realized by concrete singularities. According to one's taste, one can choose to use, say, a cusp singularity instead of two pencil singularities. What matters is that the Hopf degrees of the singularities add up to $n=\chi(M)-\tau\nu$, as is stated in \Thmref{thm-exist.sing}. One could even use just a single singularity of Hopf degree $n=\chi(M)-\tau\nu$, for example $(z_1,z_2)\maps z_1^{n+1}+z_2^2$.

In particular, this concludes the proofs of \Thmref{thm-exist.sing} and \Propref{prop-index}.

\subsection{Beyond almost-complex}
\label{sec-not.cx}

Lefschetz pencils generalize to \defemph{achiral Lefschetz pencils}. Those are Lefschetz pencils with singularities still modeled on $(z_1,z_2)\maps z_1/z_2$ and $(z_1,z_2)\maps z_1 z_2$, but this time one can use local complex coordinates that are either compatible with the orientation of $M$, or compatible with the \emph{opposite orientation} (see \cite[\S8.4]{gompf}). 

In the same spirit, the Existence \Thmref{thm-exist} can be easily generalized to a theorem that holds for more general $4$--manifolds and splittings (for cases when  $m$ from \Propref{prop-kirby} is non-zero and $(1,0)$--modifications are needed), and guarantees the existence of what we could call \defemph{achiral singular foliations}.

The surgical modification $\xi_{1,0}(q)h=qh$ can be thought of as $\Ce$--linear for the complex structure given on $\Ha$ by multiplication with complex scalars on the \emph{right} (see \Rkref{rk-cx.quaternion}(2)). The ensuing identification $\Ce^2\iso\Ha$, $(z_1,z_2)\equiv z_1+jz_2$, \emph{reverses} orientations. Nonetheless, the singularities appearing from surgical modifications with $\xi_{1,0}$ can be filled-in with the complex planes for \emph{this} complex structure, yielding a good local model (an ``anti-complex'' or ``negative'' pencil singularity). Note that any two of these complex planes will now intersect \emph{negatively}. Such a singularity can be eliminated by an anti-complex blow-up. 

More generally, of course, we can call an isolated singularity $p$ of $\F$ of \defemph{negative type} if the plane-field $\T{\F}$ on a small $3$--sphere around $p$ is preserved by a local \ac\ \str\ compatible with the \emph{opposite} orientation of $M$. Then one can define the \emph{Hopf degree} just as for singularities of positive type. The formulas from the examples above have the same degrees if we choose local complex coordinates that induce the \emph{opposite} orientation.

In conclusion, we have:

\begin{theorem}[Achiral Existence Theorem]
\label{thm-exist2}
Let $c\in H^2(M;\Zi)$ be any integral lift of $w_2(M)$, and let $c=\tau+\nu$ be any splitting. Let
\begin{align*}
& m=\rec{4}\bigl( p_1(M)+2\chi(M)-c^2 \bigr)\\
& n=\rec{4}\bigl( -p_1(M)+2\chi(M)+c^2-4\tau\nu \bigr)
\end{align*}
If $m\geq 0$ and $n\geq 0$, then 
there is an achiral singular foliation $\F$ with $e(\T{\F})=\tau$, $e(\N{\F})=\nu$, and  $m+n$ singularities. The singularities can be chosen to be modeled on the levels of the complex functions $(z_1,z_2)\maps z_1/z_2$ or $(z_1,z_2)\maps z_1 z_2$, with $n$ of them for complex-coordinates respecting the orientation, and $m$ of them for complex-coordinates reversing the orientation of $M$.

More generally, for any choice of singularities $\{p_1,\ldots,p_k\}$ of positive type and any choice of singularities $\{q_1,\ldots,q_l\}$ of negative type so that $\sum\deg p_i=n$ and $\sum\deg q_i=m$, there is an achiral singular foliation $\F$ having exactly these singularities, and with $e(\T{\F})=\tau$ and $e(\N{\F})=\nu$.
\end{theorem}

\begin{remark}
An integral lift $c$ of $w_2(M)$ is essentially a \spinc-\str. It always exists (modulo $2$--torsion in $H^2$). 
One can think of \spinc-\str s as generalizations of \ac\ \str s.
\end{remark}

\begin{example}
The $4$--sphere $\Sph{4}$ admits an achiral Lefschetz pencil with $2$--spheres meeting in a positive and a negative pencil singularity. (Fibrate $\CP^2$ by all the projective lines passing through a point, then do an anti-complex blow-down on any transverse line.)
\end{example}

The condition $m\geq 0$ was the only \emph{known} obstruction to the existence of achiral Lefschetz pencils (compare Lemma 8.4.12 in \cite{gompf}). Since achiral Lefschetz pencils are special cases of singular foliations, the theorem above adds the condition $n\geq 0$. Also, since these are the \emph{only} conditions needed for the existence of a foliation, and foliations should be expected to be much more flexible than Lefschetz pencils, this result suggests that more obstructions to the existence of achiral Lefschetz pencils probably exist and need to be uncovered.

The following obstruction to the existence of achiral singular foliations is Theorem 8.4.13 from \cite{gompf} (substituting foliations for Lefschetz pencils).

\begin{proposition}
\label{prop-achiral.obstruction}
Let $M$ be a $4$--manifold with positive-definite intersection form.
Assume that $M$ admits an achiral singular foliation $\F$ with only singular points of negative type $\{q_1,\ldots,q_k\}$. Let $m=\sum\deg q_i$. Then
\[ 1-b_1(M)+b_2(M)\geq m \]
\end{proposition}

\begin{proof}
The class $c=e(\T{\F})+e(\N{\F})$ in $H^2(M;\Zi)$ has $c\equiv w_2(M) \pmod{2}$. Thus it is a characteristic element for the intersection form: $c\cdot\alpha\equiv\alpha\cdot\alpha\pmod{2}$. 
By Donaldson's celebrated result \cite{donaldson.intersection}, a smooth $4$--manifold with positive-definite intersection form must have the intersection form 
$\bigoplus^{b_2(M)}(1)$.
Let $\{\alpha_j\}$ be any basis for the intersection form written as above. Then 
$c=\sum a_j\alpha_j +\text{torsion part}$, where
and all $a_j$ must be \emph{odd} integers. Then $c^2=\sum a_j^2\geq b_2(M)=\sigma(M)$. 
By the achiral analogue of \Propref{prop-index}, $m$ must satisfy 
$m =  \rec{4}\bigl(p_1(M)+2\chi(M)-c^2\bigr)$. Since $p_1(M)=3\sigma(M)$, we deduce that $\sigma(M)+\chi(M)\geq 2m$, which is the same as $1-b_1(M)+b_2(M)\geq m$.
\end{proof}

\begin{example}
The manifold $\#^k \,\Sph{3}\times\Sph{1}$ admits no achiral singular foliations if $k>1$. Indeed, $1-b_1+b_2=1-k<0$. (For $k=1$, we have a fibration by tori, products with $\Sph{1}$ of the circle-fibers of the Hopf fibration of $\Sph{3}$.)
\end{example}

\subsection{Closed leaves and transversals}

The strategy for proving Theorems \ref{thm-leaf} and \ref{thm-transversal} (on prescribing closed leaves and closed transversals) is the same as for the Existence \Thmref{thm-exist}: Before using Thurston's \Thmref{thm-thurston} to homotop the plane-field $\til{L}_\tau$ to a foliation, we arrange it so that it is already integrable in a certain region where either it is transversal to a certain surface or tangent to a certain surface. By keeping that region fixed, we end up with a foliation that is either transversal to the surface or has it as a leaf. 

A few small steps are necessary:

\begin{piece}
For any embedded closed surface $S$ in $M$, denote by $\nu_{S}\co\N{S}\to S$ the projection of the normal bundle of $S$. Embed $\N{S}$ as a tubular neighborhood of $S$ in $M$: $S\subset \N{S}\subset M$. 
One can pull back the bundle $\N{S}\to S$ over $\N{S}$ using $\nu_S\co\N{S}\to S$. The resulting bundle $\nu_S^*\N{S}\to\N{S}$ can then be identified with the tangent bundle to the fibers of $\N{S}$ (the vertical distribution), and thus $\nu_S^*\N{S}\subset\T{\N{S}}=\T{M}\rest{\N{S}}$. 
One can also pull $\T{S}\to S$ back over $\N{S}$ using $\nu_S$. The resulting bundle $\nu_S^*\T{S}\to\N{S}$ can be identified with a complement to $\nu_S^*\N{S}$ in $\T{\N{S}}$ (a horizontal distribution), and thus is also a subbundle of $\T{M}\rest{\N{S}}$. We thus have:  
\[ \T{M}\rest{\N{S}}=\nu_S^*\T{S}\oplus\nu_S^*\N{S} \]
\end{piece}

\begin{figure}[t]
\begin{center} 
\begin{picture}(300,110)

\put(0,50){\line(1,0){300}}
\thicklines
\put(100,51){\line(1,0){100}}

\thinlines
\multiput(105,0)(5,0){19}{\line(0,1){100}}
\put(100,0){\line(1,1){100}}

\put(150,50.3){\makebox(0,0){$\bullet$}}

\put(5,54){$M$}

\put(104.8,90){\vector(1,0){0}}
\qbezier(102,90)(80,90)(70,95)
\put(40,95){$\nu_S^*\N{S}$}

\put(101.8,2.2){\vector(1,-1){0}}
\qbezier(99,5)(80,20)(60,20)
\put(51,17.5){$s$}

\put(198,51.5){\vector(0,-1){0}}
\qbezier(197.8,54)(198,85)(225,85)
\put(227,82){$\N{S}$}

\put(152,49){\vector(-1,1){0}}
\qbezier(154,46.8)(175,25)(230,25)
\put(232,21){$S$}

\end{picture}
\end{center}
\caption{Building a complex-line bundle with $c_1=S$}
\label{fig-line.bundle}
\end{figure}
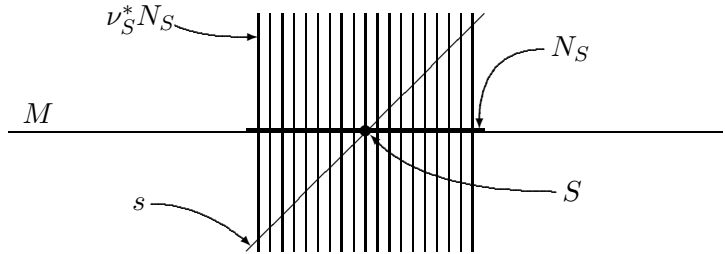

\begin{piece}
\label{pc-line.bundle}
For any surface $\Sigma$, one can built a complex-line bundle $L_\Sigma$ with Chern class $c_1(L_\Sigma)=\Sigma$ as follows: Take $\N{\Sigma}$ and pull-it back over itself using $\nu_\Sigma$. The resulting bundle $\nu_\Sigma^*\N{\Sigma}\to\N{\Sigma}$ is trivial \emph{off} $\Sigma$, since the section $s\co\N{\Sigma}\to\nu_\Sigma^*\N{\Sigma}$, $s(v)=(v,v)$ (think $\nu_\Sigma^*\N{\Sigma}\subset\N{\Sigma}\times\N{\Sigma}$), is non-zero off $\Sigma$ and hence trivializes (see \Figref{fig-line.bundle}). Therefore one can extend the bundle $\nu_\Sigma^*\N{\Sigma}$ from over $\N{\Sigma}\subset M$ to over the whole $M$, gluing it to some trivial bundle over $M\setminus \Sigma$. The result is a complex-line bundle (= oriented real $2$--plane bundle) $L_\Sigma$ with a section (an extension of $s$) that is zero only over $\Sigma$, so $c_1(L_\Sigma)=\Sigma$. Notice that, while $\nu_\Sigma^*\N{\Sigma}$ can be considered as a subbundle of $\T{M}$, in general the same is no longer true of $L_\Sigma$.
\end{piece}

\begin{figure}[t] 
\begin{center}
\begin{picture}(150,150)(0,0)

\linethickness{1.5pt}
\put(0,75){\line(1,0){150}}
\put(75,150){\line(0,-1){150}}
\thinlines
\multiput(2,100)(0,-5){11}{\line(1,0){146}}
\multiput(50,147.5)(0,-5){30}{\line(1,0){50}}

\thicklines
\put(17,130){\makebox(0,0){$\nu_A^*\N{A}$}}
\put(33,130){\vector(1,0){17}}
\put(135,130){\makebox(0,0){$\nu_S^*\T{S}$}}
\put(135,123){\vector(0,-1){23}}
\put(17,20){\makebox(0,0){$A$}}
\put(23,20){\vector(1,0){52}}
\put(135,20){\makebox(0,0){$S$}}
\put(135,27){\vector(0,1){48}}

\end{picture}
\qquad
\begin{picture}(150,150)(0,0)

\linethickness{1.5pt}
\put(0,75){\line(1,0){152}}
\put(75,153.5){\line(0,-1){152}}
\thinlines
\multiput(2,50)(7,0){15}{\line(1,1){50}}
\put(2,57){\line(1,1){43}}
\put(2,64){\line(1,1){36}}
\put(2,71){\line(1,1){29}}
\put(2,78){\line(1,1){22}}
\put(2,85){\line(1,1){15}}
\put(2,92){\line(1,1){8}}

\put(107,50){\line(1,1){43}}
\put(114,50){\line(1,1){36}}
\put(121,50){\line(1,1){29}}
\put(128,50){\line(1,1){22}}
\put(135,50){\line(1,1){15}}
\put(142,50){\line(1,1){8}}

\multiput(50,3.5)(0,7){15}{\line(1,1){50}}
\put(57,3.5){\line(1,1){43}}
\put(64,3.5){\line(1,1){36}}
\put(71,3.5){\line(1,1){29}}
\put(78,3.5){\line(1,1){22}}
\put(85,3.5){\line(1,1){15}}
\put(92,3.5){\line(1,1){8}}

\put(50,108.5){\line(1,1){43}}
\put(50,115.5){\line(1,1){36}}
\put(50,122.5){\line(1,1){29}}
\put(50,129.5){\line(1,1){22}}
\put(50,136.5){\line(1,1){15}}
\put(50,143.5){\line(1,1){8}}

\thicklines
\put(17,126){\makebox(0,0){$\nu_B^*\N{B}$}}
\put(33,126){\vector(1,0){17}}
\put(132,130){\makebox(0,0){$\nu_S^*\N{S}$}}
\put(132,123){\vector(0,-1){23}}
\put(17,20){\makebox(0,0){$B$}}
\put(23,20){\vector(1,0){52}}
\put(132,20){\makebox(0,0){$S$}}
\put(132,27){\vector(0,1){48}}

\end{picture}
\end{center}

\caption{Identifying fibers, for \ref{pc-tangent} and \ref{pc-normal}}
\label{fig-bundles}
\end{figure}
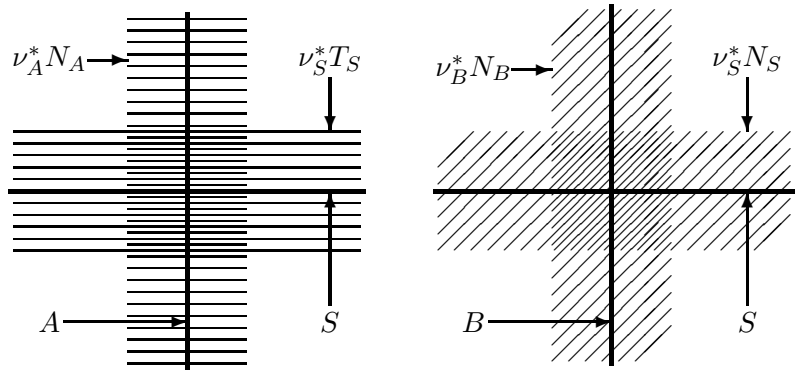

\begin{piece}
\label{pc-tangent}
Consider now an embedded \emph{connected} surface $S$ and a homology $2$--class $\alpha$ \st\ 
\[ \chi(S)=\alpha\cdot S \]
Represent $\alpha$ by an embedded surface $A$ transverse to $S$, and build a complex-line bundle $L_A$ with $c_1(L_A)=A$ as above, in \ref{pc-line.bundle}: extend $\nu_A^*\N{A}\to\N{A}$ over the whole $M$. Near $A$, the bundle $L_A$ is a subbundle of $\T{M}$ and has a section $s(v)=(v,v)$ as above. 
On the other hand, build the bundle $\nu_S^*\T{S}\to\N{S}$ as a subbundle of $\T{M}$ near $S$. Arrange that over the intersection $\N{S}\cap\N{A}$ the fibers of $\nu_S^*\T{S}$ and $\nu_A^*\N{A}$ (\ie $L_A$) coincide
(see \Figref{fig-bundles}, left). 
(Do that such that $\nu_S^*\T{S}$ is still complementary to $\nu_S^*\N{S}$ in $\T{M}\rest{\N{S}}$.) Then the section $s$ of $\nu_A^*\N{A}$ is also a section of $\nu_S^*\T{S}$ defined near $A$. Viewed there, it looks like a vector field tangent to $S$, defined only near $A$ and with zeros along $A$. But $\chi(S)=A\cdot S$, and so the zeros of $S$ along $A$ are the only obstructions to a non-zero extension of $s$ to the whole $\nu_S^*\T{S}$. We end up with a subbundle $\nu_S^*\T{S}\cup\nu_A^*\N{A}$ of $\T{M}$ over $\N{S}\cup\N{A}$, with a section $s$ that is zero only along $A$.
If we glue it to a trivial bundle over the rest of $M$, the result will be $L_A$. The difference is that now $L_A$ is a subbundle of $\T{M}$ near $S$, and is complementary there to $\nu_S^*\N{S}$. With a bit of care, we can actually get
\[ L_A\rest{\N{S}}=\nu_S^*\T{S} \]
\end{piece}

\begin{piece}
\label{pc-normal}
Consider an embedded \emph{connected} surface $S$ and a class $\beta$ \st\ 
\[ S\cdot S=\beta\cdot S \]
Then represent $\beta$ by an embedded surface $B$ transverse to $S$ and build $L_B$ as before, in \ref{pc-line.bundle}. Build also $\nu_S^*\N{S}\to\N{S}$, and arrange so that the fibers of $\nu_B^*\N{B}$ coincide to the fibers of $\nu_S^*\N{S}$ over $\N{S}\cap\N{B}$ (see \Figref{fig-bundles}, right). 
(Do that so that $\nu_S^*\N{S}$ stays complementary to $\nu_S^*\T{S}$ in $\T{M}\rest{\N{S}}$.) The section $s$ of $\nu_B^*\N{B}$ is now also a section of $\nu_S^*\N{S}$ defined near $B$. It looks like a normal vector field to $S$, defined near $B$ and with zeros on $B$. But $S\cdot S=B\cdot S$, and thus $s$ can be extended to a global section of $\nu_S^*\N{S}$ with zeros only along $B$. We end up with $\nu_S^*\N{S}\cup\nu_B^*\N{B}$, subbundle of $\T{M}$ over $\N{S}\cup\N{B}$, that is trivialized off $B$. It can be extended trivially to the whole $M$, yielding $L_B$. 
But now $L_B$ is a subbundle of $\T{M}$ near $S$, and is complementary there to $\nu_S^*\T{S}$. With a bit of care, we even get
\[ L_B\rest{\N{S}}=\nu_S^*\N{S} \]
\end{piece}

\begin{piece}
\label{pc-final}
Combine \ref{pc-tangent} and \ref{pc-normal}: Let $S$ be a \emph{connected} surface, and let $\alpha$ and $\beta$ be such that 
\[ \chi(S)=\alpha\cdot S \qquad\text{and}\qquad 
   S\cdot S=\beta\cdot S \]
Represent $\alpha$ and $\beta$ by transverse surfaces $A$ and $B$, then simultaneously re-build the bundles $L_A$ and $L_B$ in such a manner that they are complementary subbundles of $\T{M}\rest{\text{near }S}$. Namely, 
\[ L_A\rest{\N{S}}=\nu_S^*\T{S} \qquad\text{and}\qquad
  L_B\rest{\N{S}}=\nu_S^*\N{S} \]
Notice that $\nu_S^*\N{S}$ is arranged to be tangent to the fibers of $\N{S}$. Also, if $S$ has trivial normal bundle, then $\nu_S^*\T{S}$ can be arranged to be tangent to parallel copies of $S$  (and to $S$ itself) in $\N{S}$.

Notice that, if $(L_A\oplus L_B)_{m,n}\iso\T{M}$, then the resulting singular plane fields $\til{L}_A$ and $\til{L}_B$ can be kept fixed near $S$, so that $\til{L}_A\rest{\N{S}}=\nu_S^*\T{S}$ and $\til{L}_B\rest{\N{S}}=\nu_S^*\N{S}$. In general, though, they \emph{cannot} be kept fixed near $A$ or $B$. Indeed, there they must pass through an isomorphism of the type $\N{\Sigma}\oplus\T{\Sigma}\iso\til{N}\Sub{\Sigma}\oplus\uaR^2$
(since $\T{\Sigma}\oplus\uaR=\uaR^3$), where $\til{N}_{\Sigma}$ is an isomorphic copy of $\N{\Sigma}$, but a copy that is \emph{not} normal to $\Sigma$ when embedded in $\T{M}$. This moved copy $\til{N}_{\Sigma}$ becomes part of $\til{L}_\Sigma$, while $\uaR^2$ becomes part of the complementary  bundle.
\end{piece}

Finally, we are ready to assemble the above steps into the proofs of \ref{thm-transversal} and \ref{thm-leaf}:

\begin{Proof}{Transversal \Thmref{thm-transversal}}
The statement we need to prove is:

\smallskip
\emph{Let $S$ be a closed connected surface. 
Let $c$ be a complex class with a splitting  $c=\tau+\nu$ \st\ $\chi(M)-\tau\nu\geq 0$. 
If $\chi(S)=\nu\cdot S$ and $S\cdot S=\tau\cdot S$
then there is a singular foliation $\F$ with $e(\T{\F})=\tau$, $e(\N{\F})=\nu$, and having $S$ as a closed transversal.
}

\smallskip
Build the line bundles $L_\tau$ and $L_\nu$ following the recipe from \ref{pc-final}. Do the surgical modifications on $L_\tau\oplus L_\nu$ far from $S$. The resulting singular plane fields $\til{L}_\tau$ and $\til{L}_\nu$ now have $\til{L}_\tau\rest{\N{S}}=\nu_S^*\N{S}$, and thus $\til{L}_\tau$ can be arranged to be tangent to the fibers of $\N{S}$ in $M$ (in other words, $\til{L}_\tau$ is integrable near $S$).
Keeping the plane field $\til{L}_\tau$ fixed near the filled-in singularities and near $S$, we end up, after applying Thurston's \Thmref{thm-thurston}, with a singular foliation $\F$ having the fibers of $\N{S}$ as pieces of leaves. Thus $S$ is everywhere transverse to the $\F$.
\end{Proof}

\begin{Proof}{Leaf \Thmref{thm-leaf}}
The statement we need to prove is:

\smallskip
\emph{Let $S$ be a closed connected surface with $S\cdot S\geq 0$. 
Let $c$ be a complex class with a splitting $c=\tau+\nu$ \st\  $\chi(M)-\tau\nu\geq S\cdot S$.
If $\chi(S)=\tau\cdot S$ and $S\cdot S=\nu\cdot S$
then there is a singular foliation $\F$ with $e(\T{\F})=\tau$, $e(\N{\F})=\nu$, and having $S$ as a closed leaf.
(The number of singularities along $S$ is $S\cdot S$.)
}

\smallskip
{\bf A.}
Assume first that the normal bundle $\N{S}$ of $S$ is trivial.
Build the line bundles $L_\tau$ and $L_\nu$ following the recipe from \ref{pc-final}. Do the surgical modifications on $L_\tau\oplus L_\nu$ far from $S$. The resulting singular plane fields $\til{L}_\tau$ and $\til{L}_\nu$ now have $\til{L}_\tau\rest{\N{S}}=\nu_S^*\T{S}$, and, since $\N{S}$ is trivial, $\til{L}_\tau$ can be arranged to be tangent to parallel copies of $S$. Keeping the plane field $\til{L}_\tau$ fixed near the filled-in singularities and near $S$, we end up, after applying Thurston's \Thmref{thm-thurston}, with a singular foliation with $S$ (and its parallel copies) as leaves.

\begin{figure}[t]
\begin{center} 
\begin{picture}(300,70)(-150,-35)

\thicklines
\put(-150,0){\line(1,0){300}}

\thinlines

\put(-150,15){\line(1,0){100}}
\put(150,-15){\line(-1,0){100}}

\qbezier(-50,15)(-22,15)(-15,10)
\qbezier(50,-15)(22,-15)(15,-10)

\put(-150,30){\line(1,0){120}}
\put(150,-30){\line(-1,0){120}}

\qbezier(-30,30)(-20,30)(-10,15)
\qbezier(30,-30)(20,-30)(10,-15)
\
\put(0,0){\line(0,1){35}}
\put(0,0){\line(0,-1){35}}

\put(0,0){\line(-2,3){10}}
\put(0,0){\line(2,-3){10}}
\put(0,0){\line(-2,-3){20}}
\put(0,0){\line(2,3){20}}

\put(0,0){\line(-3,2){15}}
\put(0,0){\line(-3,-2){30}}
\put(0,0){\line(3,-2){15}}
\put(0,0){\line(3,2){30}}

\put(130,3){$S$}

\end{picture}
\end{center}
\caption{Foliating around $S$ when $\N{S}$ is non-trivial}
\label{fig-singleaf}
\end{figure}
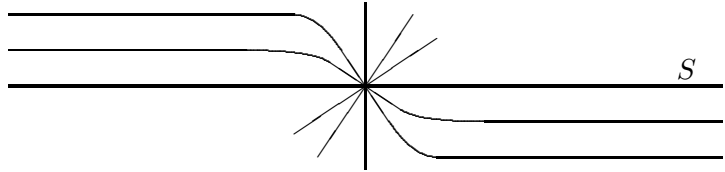

\smallskip
{\bf B.}
In general, if $\N{S}$ is not trivial, then we will place pencil singularities along $S$, as suggested in \Figref{fig-singleaf}. 
Having foliated a neighborhood of $S$, we can essentially apply the same recipe as above.

Notice that the condition $\chi(M)-\tau\nu\geq S\cdot S$ is there merely to ensure that we have enough singularities available. 
We leave the remaining details of the proof of the Leaf \Thmref{thm-leaf} to the elusive interested reader.
\end{Proof}

If one starts with a surface with $S\cdot S<0$, then one could try to use negative singularities to foliate a neighborhood.  Thus, one needs achiral foliations:

\begin{proposition}
\label{prop-achiral.leaf}
Let $S$ be a closed connected surface with $S\cdot S<0$. Let $c, \tau, \nu, m, n$ be as in \ref{thm-exist2}. If $m\geq -S\cdot S$ and $n\geq 0$, then there is an achiral singular foliation having $S$ as a leaf.
\end{proposition}

\section{Appendix}
\label{sec-appendix}

\begin{Proof}{\Lemmaref{lemma-domega}}
We prove that, if $g$ is a \riem\ metric, $\nabla$ 
its \LC\ connection, $J$ be any $g$--orthogonal \ac\ \str, and
$\omega(x,y)=\inner{Jx,\,y}$ its fundamental $2$--form,
then, for any vector fields $x, z$ on $M$, we have:
\[ (d\omega)(x,Jx,z)
 =\Inner{[x,Jx],\,Jz}
 -\Inner{\Nabla{x}x+\Nabla{Jx}Jx,\ z} \]

For any $2$--form $\alpha$ we have:
\begin{align*}
(d\alpha)(x,y,z)
&=\bigl(\Nabla{x}\alpha\bigr)(y,z)
 +\bigl(\Nabla{y}\alpha\bigr)(z,x)
 +\bigl(\Nabla{z}\alpha\bigr)(x,y)\\
&=x\alpha(y,z)+y\alpha(z,x)+z\alpha(x,y)\\
&\qquad
 -\alpha\bigl(\Nabla{x}y,\,z\bigr)-\alpha\bigl(\Nabla{y}z,\,x\bigr)
 -\alpha\bigl(\Nabla{z}x,\,y\bigr)\\
&\qquad
 -\alpha\bigl(y,\,\Nabla{x}z\bigr)-\alpha\bigl(z,\,\Nabla{y}x\bigr)
 -\alpha\bigl(x,\,\Nabla{z}y\bigr)
\end{align*}
Applying this to $\omega(a,b)=\inner{Ja,\,b}$, we have:
\begin{align*}
(d\omega)(x,Jx,z)
&=-x\Inner{x,z}+(Jx)\Inner{Jz,\,x}+z\Inner{Jx,Jx}\\
&\qquad
 -\Inner{J\Nabla{x}Jx,\,z}-\Inner{J\Nabla{Jx}z,\,x}
 -\Inner{J\Nabla{z}x,\,Jx}\\
&\qquad
 +\Inner{x,\,\Nabla{x}z}-\Inner{Jz,\,\Nabla{Jx}x}
 -\Inner{Jx,\,\Nabla{z}Jx}
\end{align*}
Using that $\inner{a,Jb}=-\inner{Ja,b}$ we get: 
\begin{align*}
(d\omega)(x,Jx,z)
&=-x\Inner{x,z}+(Jx)\Inner{Jz,\,x}+z\Inner{x,x}\\
&\qquad
 +\Inner{\Nabla{x}Jx,\,Jz}+\Inner{\Nabla{Jx}z,\,Jx}
 -\Inner{\Nabla{z}x,\,x}\\
&\qquad
 +\Inner{x,\,\Nabla{x}z}-\Inner{Jz,\,\Nabla{Jx}x}
 -\Inner{Jx,\,\Nabla{z}Jx}
\end{align*}
Since 
$\Inner{\Nabla{z}x,\,x}=\rec{2}z\Inner{x,x}$
and
$\Inner{Jx,\,\Nabla{z}Jx}=\rec{2}z\Inner{Jx,Jx}=\rec{2}z\Inner{x,x}$, 
we cancel the last terms of  each line, and get:
\begin{align*}
(d\omega)(x,Jx,z)
&=-x\Inner{x,z}+(Jx)\Inner{Jz,\,x}\\
&\qquad
 +\Inner{\Nabla{x}Jx,\,Jz}+\Inner{\Nabla{Jx}z,\,Jx}\\
&\qquad
 +\Inner{x,\,\Nabla{x}z}-\Inner{Jz,\,\Nabla{Jx}x}
\end{align*}
But $x\Inner{x,z}=\Inner{\Nabla{x}x,\,z}+\Inner{x,\,\Nabla{x}z}$, so
$(Jx)\Inner{Jz,\,x} = -(Jx)\Inner{z,\,Jx}
 =
 \break 
 -\Inner{\Nabla{Jx}z,\,Jx}-\Inner{z,\,\Nabla{Jx}Jx}$, 
and therefore:
\begin{align*}
(d\omega)(x,Jx,z)
&=-\Inner{\Nabla{x}x,\,z}-\Inner{x,\,\Nabla{x}z}
 -\Inner{\Nabla{Jx}z,\,Jx}-\Inner{z,\,\Nabla{Jx}Jx}\\
&\quad
 +\Inner{\Nabla{x}Jx,\,Jz}+\Inner{\Nabla{Jx}z,\,Jx}
 +\Inner{x,\,\Nabla{x}z}-\Inner{Jz,\,\Nabla{Jx}x}\\
(d\omega)(x,Jx,z)
&=-\Inner{\Nabla{x}x,\,z}-\Inner{z,\,\Nabla{Jx}Jx}
 +\Inner{\Nabla{x}Jx,\,Jz}-\Inner{Jz,\,\Nabla{Jx}x}
\end{align*}
Since $\nabla$ is torsion-free, we have 
$\Nabla{x}Jx-\Nabla{Jx}x=[x,Jx]$, so:
\[ (d\omega)(x,Jx,z)
 =\Inner{[x,Jx],\,Jz}
 -\Inner{\Nabla{x}x+\Nabla{Jx}Jx,\ z} \]
which concludes the proof.
\end{Proof}%

In particular, if $\omega$ is symplectic (\ie $d\omega=0$), then any 
integrable $J$--holomor\-phic plane field is $g$--minimal, and, vice-versa, any $g$--minimal 
$J$--holomorphic plane field must be integrable. 
The converse is also true: If there are enough $J$--holomorphic integrable minimal 
plane fields, then $\omega$ must be symplectic. 
Thus:

\begin{corollary}
Assume that $M$ admits two transversal $2$--dimensional foliations 
$\F$ and $\G$ such that: 
there is a metric $g$ \st\ both $\F$ and $\G$ are 
$g$--minimal,
and there is a $g$--orthogonal \ac\ \str\ $J$ that makes both $\F$ 
and $\G$ be $J$--holomorphic.
Then $M$ admits the symplectic structure 
$\omega(x,y)=g(Jx,y)$.
\end{corollary}

In particular, if the first condition is satisfied, and further 
$\F$ and $\G$ are $g$--orthogonal, then the second condition is 
automatically satisfied, and \Propref{prop-2taut.sympl} follows:

{\sl
If a \riem\ manifold $M$ admits two $g$--orthogonal and $g$--minimal 
foliations, then $M$ admits a symplectic structure.
}


\Addresses\recd

\end{document}